\documentclass[12pt]{amsart}
\usepackage{graphicx, overpic}
\usepackage[below]{placeins}
\DeclareMathOperator{\sech}{sech}

\usepackage{amsmath,amscd,amssymb, epsfig,psfrag}

\newtheorem{theorem}{Theorem}[section]

\newtheorem{lem}[theorem]{Lemma}
\newtheorem{prop}[theorem]{Proposition}
\newtheorem{cor}[theorem]{Corollary}
\theoremstyle{definition}
\newtheorem{definition}[theorem]{Definition}

\newtheorem{conjecture}[theorem]{Conjecture}

\newcommand{\NN}{\mathbb{N}}

\newcommand{\CC}{\mathbb{C}}
\newcommand{\RR}{\mathbb{R}}

\def\vol{\mbox{\rm{Vol}}}

\def\a{\alpha}

\def\D{\Delta}

\def\g{\gamma}

\def\l{\lambda}
\def\k{\kappa}

\def\S{\Sigma}
\def\dd{\delta}
\def\eg{{\it e.g.}}

\newcommand{\U}{\ensuremath{\widetilde}}

\def\Vol{\mbox{\rm{Vol}}}
\def\min{\mbox{\text{min}}}

\def\HH{\mathbb{H}}

\def\MM{\mathcal{M}}
\def\PP{\mathbb{P}}
\def\CC{\mathbb{C}}

\def\split{\backslash\backslash}

\edef\t@mp{\catcode`\noexpand\#=\the\catcode`\#}%
    \catcode`\#=12
    \def\h@sh{#}%
    \t@mp

\edef\t@mp{\catcode`\noexpand\~=\the\catcode`\~}%
    \catcode`\~=12
    \def\tild@{~}%
    \t@mp

\begin{document}

\title{Lower bounds on volumes of hyperbolic Haken 3-manifolds} 

\author[Ian Agol]{%
        Ian Agol} 
\address{%
        MSCS UIC 322 SEO, m/c 249\\
        851 S. Morgan St.\\
        Chicago, IL 60607-7045} 
\email{%
        agol@math.uic.edu}  
\thanks{Agol partially supported by NSF grant DMS-0204142 and the Sloan
Foundation}

\author[Pete Storm]{%
	Peter A. Storm}
\address{%
Stanford University\\
Mathematics, Bldg. 380\\
450 Serra Mall\\
Stanford, CA 94305-2125
}

\email{%
storm@math.stanford.edu
}

\thanks{Storm partially supported by an NSF postdoctoral fellowship}

\author[Bill Thurston]{%
	William P. Thurston\\ with an appendix by Nathan Dunfield}
\address{%
310 Malott Hall \\
Cornell University\\
Ithaca, NY 14853-4201 
}
\email{%
wpt@math.cornell.edu
}

\thanks{Thurston partially supported by the NSF grant DMS-0343694}
\address{%
Mathematics 253-37
Caltech
Pasadena, CA 91125
}
\email{%
dunfield@caltech.edu}

\thanks{Dunfield partially supported by the NSF grant DMS-0405491 and the Sloan foundation}

\date{%
 June 16, 2005}


\begin{abstract} 
We prove a volume inequality for 3-manifolds having $C^{0}$ metrics ``bent'' along
a hypersurface, and satisfying certain curvature pinching conditions. The result
makes use of Perelman's work on Ricci flow and geometrization of closed 3-manifolds.
Corollaries include a new proof of a conjecture of Bonahon about volumes of 
convex cores of Kleinian groups, improved volume estimates for certain Haken hyperbolic
3-manifolds, and a lower bound on the minimal volume orientable hyperbolic 3-manifold. 
An appendix by Dunfield compares estimates of volumes of hyperbolic
3-manifolds drilled along a closed embedded geodesic with experimental data. 
\end{abstract} 

\maketitle
\section{Introduction}

In this paper, we prove a volume inequality for 3-manifolds having $C^{0}$ metrics ``bent'' 
along a surface 
satisfying certain curvature pinching conditions. The method is a direct generalization
of work of Bray and Miao on the Penrose conjecture. We also make use of Perelman's results
on geometrization \cite{Per02,Per03}. Perelman's montonicity formula for the Ricci flow with surgery 
(discovered by Hamilton in the case of standard Ricci flow) implies that if 
 $(M,g)$ is a hyperbolic 3-manifold, and $(M,h)$ is a Riemannian metric such
that the scalar curvature $R(h)\geq -6=R(g)$, then $\vol(M,h)\geq \vol(M,g)$. 
More generally, if $M$ is not hyperbolic, then $\vol(M,h)\geq V_{3}\|M\|$, where $V_{3}\|M\|$ 
denotes the simplicial volume of $M$ ($V_{3}$ is the volume of a regular ideal tetrahedron
in $\HH^{3}$, and $\|M\|$ is the Gromov norm of the fundamental class of the 3-manifold $M$).
This sort of pinching condition is much weaker
than pinching conditions using Ricci curvature coming from the work of Besson, Courtois,
and Gallot, which implies a similar volume estimate if $Ric(h)\geq -2 h$. 

Theorem \ref{main inequality} is the main result  of this paper, which states that if $(M,g)$ is a compact hyperbolic 3-manifold with minimal
surface boundary, then $\vol(M,g)\geq \frac12 V_{3}\|DM\|$, where $DM$ denotes
the double of $M$ along its boundary. There are some interesting consequences of
this result. First, it implies a strong form of a conjecture of Bonahon, which states
that the volume of a hyperbolic 3-manifold $M$ with convex boundary 
is $\geq \frac12 V_{3}\|DM\|$, with equality only in the case that $M$ has geodesic
boundary. This was proven by the second author in \cite{Storm04}, using the methods of Besson, Courtois, and Gallot \cite{BCG}. 
Since one may find a minimal surface representing 
the (maximally compressed) boundary of $M$, Bonahon's conjecture follows immediately from 
theorem \ref{main inequality} (see theorem \ref{noncompactbonahon}). 

Theorem \ref{main inequality} also implies that if $M$ is a compact hyperbolic Haken 3-manifold (without
boundary), and $\S\subset M$ is an incompressible surface, then 
$\vol(M)\geq \frac12 V_{3}\|D(M-\mathcal{N}(\S))\|$ (where $\mathcal{N}(\S)$ denotes a regular
neighborhood of $\S$). This follows immediately from theorem \ref{main inequality}
by noticing that $\S$ is isotopic to a minimal surface in $M$. We also generalize this result to
the finite volume case in theorem \ref{noncompact}. One illuminating special case is to consider
a compact hyperbolic 3-manifold $M$ with geodesic boundary. If one doubles $M$, there is a
natural hyperbolic structure on $DM$ induced from $M$. If one takes two copies of $M$, and
glues the two copies of $\partial M$ by a diffeomorphism $\psi $ which is not isotopic to an
isometry of $\partial M$, then by Thurston's hyperbolization theorem  $M\cup_{\psi} M$ has a canonical
hyperbolic metric. Theorem \ref{main inequality} implies that $\vol(M\cup_{\psi} M) > \vol(DM)$. 
One may compare this to gluing manifolds along tori, in which case the simplicial volume is additive \cite{Th}.

Another application of our volume inequality given in section \ref{minvol} is to show that the minimal volume closed orientable hyperbolic
3-manifold has volume $\geq .67$ (due to Agol and Dunfield). This improves on the previous best lower
bound of $.33$, due to Przeworski \cite{Pr3}. It is conjectured that the Weeks manifold,
with volume $.9427...$ is the minimal volume orientable hyperbolic 3-manifold. If $M$ is a closed hyperbolic 3-manifold  with a closed embedded
geodesic $\g$ with tube radius $R$, then we estimate the hyperbolic volume of $M-\g$ in terms of $\vol(M)$ and $R$. This improves on the main 
result of \cite{Ag0}. In the appendix, Dunfield compares these estimates with experimental data
via Goodman's program {\it tube} \cite{Snap}.

In section 2, we state definitions and the main theorem, as well as some immediate corollaries.
In section 3, we state Perelman's monotonicity results. Sections 4-7 prove the main theorem,
and sections 8 and 9 extend this to the non-compact case. Section 10 proves the results on minimal
volume orientable hyperbolic 3-manifolds, and section 11 discusses universal manifold pairings. Section
12 gives some open questions stemming from this work. 

{\bf Acknowledgements:}
We thank Pengzi Miao, Andres Neves, and Rick Schoen
for  helpful correspondence.

\section{Definitions and statement of the main theorem}
Usually, we will be assuming that manifolds are orientable in this paper (most results
for the non-orientable case follow by passing to the 2-fold orientable cover). 
A properly embedded 
{\it incompressible} surface $S$ in $M^3$ is a
surface for which the fundamental group injects (excluding
$S^2$ and $\RR\PP^{2})$. A manifold
is {\it irreducible} if every embedded $S^2$ bounds a ball. An irreducible
manifold with an incompressible surface is called {\it Haken}.  
If a Riemannian manifold has boundary, then it is hyperbolic with totally
geodesic boundary if the metric is locally modelled on a closed half-space 
in $\HH^{3}$ bounded by a geodesic plane. If $(M^{3},g)$ is a Riemannian manifold,
and $\S^{2}\subset M$ is an embedded surface, then $M\split \S$ denotes
the Riemannian manifold with boundary obtained by taking the path metric
completion of $M-\S$. $M\split \S$ will have new boundary corresponding
to the unit normal bundle of $\S\subset M$. 

The following theorem is a special case of theorem \ref{main inequality}:

{\bf Theorem }
{\it If $(M,g)$ is a compact hyperbolic 3-manifold with minimal surface boundary, 
then $\vol(M,g)\geq  \frac12 V_{3}\|DM\|$. }

In particular, if $M$ is acylindrical, then $M$ admits
a hyperbolic metric $\nu$ with totally geodesic boundary, and 
$\vol(M)\geq \vol(M,\nu)$. We conjecture that in theorem \ref{main inequality}, we do not need to assume
that the boundary is compact.

Let $V_{8}= 3.66...$ denote the volume of a regular
ideal octahedron in $\HH^{3}$. Miyamoto showed that if $M$ is hyperbolic 
with totally
geodesic boundary, then $\vol(M)\geq -V_{8}\chi(M)$, with
equality holding only for manifolds composed of regular ideal
octahedra, glued together in the pattern of an ideal
triangulation \cite{Mi}.

\begin{definition}
If $M$ is a manifold with boundary whose interior admits a hyperbolic structure, 
$guts(M)=N \subset M$ is a submanifold
such that $\partial N =\partial_{0} N \cup \partial_{1}N$, where $\partial_{i} N$ is a subsurface, 
$\partial_{0} N=N\cap \partial M$, $\partial_{1} N$ consists of annuli or tori such that 
$\partial\partial_{1} N = \partial_{1} N\cap \partial M$.  
Moreover, $(N,\partial_{1} N)$ is the maximal pared acylindrical submanifold such that
no components of $N$ are solid tori. One may also characterize $guts(M)$ by
$N-\partial_{1}N$ admits a complete hyperbolic
metric with totally geodesic boundary, and $D(M-N)$ admits a graph manifold structure \cite{Mo, Johannson79, JS79}.  
\end{definition}

The following result was conjectured in \cite{Ag3}. 

\begin{cor} \label{guts} Let $M$ be a closed hyperbolic 3-manifold, and  $\S\subset M$ be
an incompressible surface. Then
$\vol(M)\geq \frac12 V_{3}\| D(M\split \S)\|\geq  -V_{8}\chi (guts(M\split \S))$.
\end{cor}

In particular, if $(M,g)$ is a Haken hyperbolic 
3-manifold such that $\vol(M,g)\leq V_{8}$, then for any incompressible surface
$\S\subset M$, $guts(M\split \S) =\emptyset$, which means that $\S$ is a so-called fibroid
(this improves greatly on the estimates of \cite{CS1, H}). The ramifications of this corollary
will be pursued in subsequent papers.  We record for now the following improvement
on a theorem of Lackenby. If $D$ is the projection of a link $L\subset S^{3}$, then 
the crossings of $D$ are divided into {\it twist} equivalence classes, where two crossings 
are twist equivalent if there is a loop in the projection plane intersecting $D$ transversely precisely 
in the two crossings. The number of twist classes of crossings $t(D)$ is called the twist number
of $D$. One assumes that the diagram $D$ has no crossings separating $D$. 
\begin{theorem} Let $L\subset S^{3}$ be a hyperbolic alternating link with alternating projection $D$. 
Then 
$$V_{8}(t(D)/2-1) \leq \vol(S^{3}-L) \leq 10 V_{3}(t(D)-1).$$
\end{theorem}
The lower inequality is an equality for the Borromean rings, and follows directly from 
theorem 5 of \cite{Lackenby04} in conjunction with corollary \ref{guts}. The upper
inequality is given in the appendix to \cite{Lackenby04}.

\section{Monotonicity formulae for the Ricci flow} \label{monotonicity section}

If $(M,g)$ is a closed Riemannian manifold, let $\l(g)$ denote the minimal eigenvalue of the
operator $-4\D_{g}+R(g)$, and define the scale invariant quantity 
$$V_{\l}(g)=\vol(M, g) \left(-\frac16 \min\{\l(g),0\}\right)^{3/2}$$
(which is equivalent  to taking the volume of the metric rescaled so that $\l(g)=-6$ when $\l(g)<0$).   
Recall that $V_{3}\|M\|$ is the sum of the volumes of the hyperbolic pieces of the
geometric decomposition of $M$ \cite{Th}. 
If $R_{min}(g)\geq 0$ (where $R_{min}(g)$
is the minimum of $R(g)$), then $\|M\|=0$ \cite{SchoenYau79}. Also, define 
$$V_{R}(g)=\vol(M,g)\left(-\frac16 \min\{R_{min}(g),0\}\right)^{3/2}.$$
We have the inequality $\l(g)\geq R_{min}(g)$, as can be seen by neglecting the $-\D_{g}$ term
in the Rayleigh quotient, and thus $V_{R}(g) \geq V_{\l}(g)$. Perelman showed that $V_{R}(g)$ and
$V_{\l}(g)$ are monotonic decreasing for metrics evolving by the Ricci flow with surgery (monotonicity
of $V_{R}(g)$ for Ricci flow was shown by Hamilton \cite{Ham99}). 
Then Perelman's analysis of the geometric decomposition forming at $\infty$ time  under Ricci flow with surgery implies the following
theorem (a restatement of theorem 8.2 \cite{Per02}). 
\begin{theorem} \label{monotonicity}
Let $(M,g)$ be a closed 3-manifold. Then
$$V_{R}(g)\geq V_{\l}(g)\geq V_{3}\|M\|.$$
\end{theorem}
It is much simpler to prove
that $V_{R}(g)$
is monotonic decreasing with respect to Ricci flow with surgery, than to prove the  monotonicity of $V_{\l}(g)$.
Specifically, the above monotonicity formula for $V_{\l}(g)$ through surgeries has not been discussed in the notes \cite{KL03}, whereas the monotonicity
of $V_{R}(g)$ through surgery is obvious from Perelman's formulation of surgery, since the surgery decreases
volume, and does not change $\min\{R_{min}(g),0\}$, since the surgeries occur at a part of the manifold where
$R(g)\gg 0$.  Thus, in the
proof of theorem \ref{main inequality}, we give two different arguments, one estimating $\l(g)$, and the
other performing a conformal change to estimate $R_{min}(g)$. Although the authors of this paper have not
checked the validity of Perelman's papers, there are references available which work out many of the
details of Perelman's arguments (and fix some errors, see the website maintained by Kleiner \cite{KL03}). 
Also, a key aspect of Perelman's geometrization proof is a collapsing
theorem which he has not made available. Since this is a key to the above theorem \ref{monotonicity}, we 
note that Shioya and Yamaguchi have a paper \cite{ShY03} which claims to give an alternative proof of the collapsing
results needed by Perelman.

\section{Smoothing with a lower scalar curvature bound} \label{scalar section}

This section produces smooth Riemannian approximations to a particular kind of singular Riemannian metric while maintaining a lower bound on the scalar curvature.  The main result is proposition \ref{scalar prop}.  In a different setting, this smoothing procedure was developed as a tool to solve the Riemannian Penrose
conjecture in general relativity (see Bray \cite{Bray01} and Miao \cite{Miao02}).  The techniques of this section are adaptations of their methods to the present situation.

Let $M_+$ and $M_-$ be smooth Riemannian $3$-manifolds with compact boundary, equipped with a diffeomorphism $\Phi: \partial M_+ \longrightarrow \partial M_-$.  Fix a smooth Riemannian metric $g_{\pm}$  on $M_{\pm}$.  For a point $x \in \partial M_{\pm}$ let $H(\partial M_\pm, g_\pm) (x)$ denote the mean curvature of $\partial M_\pm$ with respect to an outward pointing normal vector field.  Let us assume that 
\begin{equation} \label{gluing condition}
\Phi^* g_- |_{\partial M_-} =  g_+ |_{\partial M_+} \quad \text{and} \quad H(\partial M_+, g_+) (x) + H( \partial M_-, g_-) (\Phi (x)) \ge 0
\end{equation}
for all $x \in \partial M_+$ (that is, the boundaries are identified by an isometry, and the mean curvatures pointwise sum to $\geq 0$).  Two relevant examples are: $M_{\pm}$ are two copies of a hyperbolic $3$-manifold with minimal surface boundary, or a finite volume hyperbolic $3$-manifold truncated along a tube about a geodesic together with a horospherical cusp of appropriate shape.

Glue $M_+$ to $M_-$ along their boundary via the map $\Phi$ to form
\begin{equation} \label{glue}
 (M, g) :=  (M_+, g_+) \bigcup_\Phi  (M_-, g_-).
\end{equation}
The manifold $M$ is smooth.  The metric tensors $g_+$ and $g_-$ agree along the boundary and glue together to form a piecewise smooth continuous Riemannian metric on $M$.  In particular, $g$ is a Lipschitz Riemannian metric (this means that in the $C^{\infty}$ structure on $M$, in
any chart the metric $g$ will a Lipschitz function on $\RR^{3}$).  We will use $\Sigma \subset M$ to denote the submanifold formed by the boundary of $M_+$ and $M_-$.  Fix a collared tubular neighborhood $\Sigma \times (- \delta_0, \delta_0) \hookrightarrow M$ taking $\Sigma \times \{ 0 \}$ to $\Sigma \subset M$.

Under these hypotheses, Miao has proven

\begin{prop} \cite[proposition 3.1]{Miao02} \label{lowerbound}
Suppose that $(M,g)$ is defined by equation \ref{glue}, satisfying the conditions of equation \ref{gluing condition}. There exists a family of smooth Riemannian metrics $\{ g_\delta \}_{0<\delta<\delta_0}$ on $M$ so that $g_\delta$ converges uniformly to $g$ as $\delta \rightarrow 0$, the metrics $g_\delta$ and $g$ agree outside $\Sigma \times (-\delta/2, \delta/2)$, and the scalar curvature $R(g_\delta)$ is bounded below by a constant $\mathcal{S}$ independent of $\delta$.
\end{prop}

{\bf Remark:} Miao performs the smoothing via  a {\it mollification} method. In the main case of interest in this
paper which have an extra symmetry, one may use the more elementary method of Bray \cite{Bray01}. 

{\noindent}Stated more geometrically, the fact that the metrics $g_\delta$ converge uniformly to $g$ implies that the Lipschitz constant of the identity map $(M,g_\delta) \rightarrow (M,g)$ converges to one.  In particular, the volume of $(M,g_\delta)$ converges to the volume of $(M,g)$.

The goal of this section is to prove

\begin{prop} \label{scalar prop}
Assume that $R(g_{\pm}) \geq  -6$.  Then there exists a family of smooth Riemannian metrics $G_\delta$ on $M$ such that the metrics $G_\delta$ converge uniformly to $g$ and have scalar curvature $R(G_{\dd})$ bounded below by $-6$.
\end{prop}

The metrics $G_\delta$ will be obtained by conformally altering the smooth metrics $g_\delta$ produced in proposition \ref{lowerbound}.  The method is an adaptation of techiques of Bray \cite{Bray01} and Miao \cite{Miao02} dealing with the case of open asymptotically flat manifolds of nonnegative scalar curvature.

Before beginning the proof, we note that if for all sufficiently small $\delta$,   $R(g_\delta)\geq -6$, then we are done.  So we can assume without a loss of generality that $R(g_\delta) \geq \mathcal{S}$, such that  $\mathcal{S}<-6$, where $\mathcal{S}$ is independent of $\delta$ by proposition \ref{lowerbound}.

Define a smooth function $\mathcal{R}_\delta: M \longrightarrow [\mathcal{S}, 0]$ such that $\mathcal{R}_\delta$ restricted to $\Sigma \times (-\delta/2, \delta/2)$ is exactly $\mathcal{S}$, and $\mathcal{R}_\delta$ vanishes outside $\Sigma \times (-\delta, \delta)$.  Note that the $g_\delta$-volume of the support of $\mathcal{R}_\delta$ converges to zero.

Define the operator $-8 \Delta_\delta + \mathcal{R}_\delta$, where $\Delta_\delta$ is the Laplacian for the metric $g_\delta$.  This Schr{\"o}dinger operator has a minimal eigenvalue $\lambda_\delta$ with a nonnegative smooth eigenfunction $u_\delta$.  Let us normalize the eigenfunction so that
$$\int_M u_\delta^2 \, dg_\delta = \vol(M,g_{\dd}).$$
If $\lambda_\delta$ were nonnegative then the eigenfunction would satisfy
$$-8 \Delta u_\delta = \lambda_\delta u_\delta - \mathcal{R}_\delta u_\delta \ge 0,$$
making $u_\delta$ a superharmonic function.  Since $M$ is compact and $\mathcal{R}_\delta$ is not a constant function, this is impossible.  Therefore $\lambda_\delta < 0$.
Take the equation 
$$-8 \Delta u_\delta + \mathcal{R}_\delta u_\delta = \lambda_\delta u_\delta,$$
multiply it by $u_\delta$, integrate over $M$, and 
perform integration by parts to obtain
\begin{equation} \label{lambda estimate}
\int_M (8 | \nabla u_\delta |^2 + \mathcal{R}_\delta u_\delta^2) \, dg_\delta = \lambda_\dd \vol(M,g_{\dd}) < 0.  
\end{equation}
This yields the estimate
$$\int_M | \nabla u_\delta |^2 \, dg_\delta  <  \frac{1}{8} \int_M  | \mathcal{R}_\delta | \, u_\delta^2 \, dg_\delta  
			\le \frac{| \mathcal{S}|}{8}\vol(M,g_{\dd}).$$

The following lemma is based on the method of proposition 4.1 \cite{Miao02}. 
\begin{lem} \label{Sobolev}
$\int_M \mathcal{R}_\delta u_\delta^2  \, dg_\delta \longrightarrow 0.$
\end{lem}
\begin{proof}
Note that  $| \mathcal{R}_\delta |$ is uniformly bounded by $| \mathcal{S} |$.
Apply H{\"o}lder's inequality to obtain
$$ \int_M | \mathcal{R}_\delta | \, u_\delta^2 \, dg_\delta \le
	\left( \int_M | \mathcal{R}_\delta |^{3/2} \, dg_\delta \right)^{2/3} \cdot 
	\left( \int_M u_\delta^6 \, dg_\delta \right)^{1/3}\leq |\mathcal{S}| (\vol(supp(\mathcal{R}_{\dd}), g_{\dd}))^{2/3} \| u_{\dd}\|_{6}^{2}. $$
If we can show that $\|u_{\dd}\|_{6}$ is uniformly bounded, then the term to the right $\to 0$, 
since the volume of the support of $\mathcal{R}_\delta$ goes to zero. 
We now apply the Sobolev inequality (see \eg \cite[Thm. 2.21]{Aubin98})  to obtain
$$\|u_{\dd}\|_{6} \le C \left( \| \nabla u_{\dd}\|_{2}+ \|u_{\dd}\|_{2} \right) < C\, \vol(M,g_{\dd})^{\frac12} (  (| \mathcal{S} |/8)^{\frac12} + 1)$$
Note that the uniform convergence of the metrics $g_\delta$ to $g$ implies we can use a single Sobolev constant $C$ for all the metrics $g_\delta$.  Therefore the number on the right hand side of the Sobolev inequality is uniformly bounded for all $\delta$, which completes the proof that $\|u_{\dd}\|_{6}$
is bounded. 
\end{proof}

This lemma has two corollaries.

\begin{cor}
$$\int_M | \nabla u_\delta |^2 \, dg_\delta \longrightarrow 0 \quad \text{and} \quad \lambda_\delta \rightarrow 0.$$
\end{cor}
\begin{proof}
Since $| \nabla u_\delta |^2 \ge 0$, the corollary follows from the lemma and equation \ref{lambda estimate}.
\end{proof}

Define the average
$$\overline{u}_\delta := \frac{1}{\text{Vol}(M, g_\delta)} \int_M u_\delta \, dg_\delta.$$

\begin{cor}
$$\int (u_{\dd}-\overline{u}_{\dd})^{2} dg_{\dd} \to 0,\quad  \overline{u}_{\dd} \to 1.$$
\end{cor}
\begin{proof}
The Poincar{\'e} inequality for $(M,g_\delta)$ yields
$$\int_M ( u_\delta - \overline{u}_\delta )^2 \, dg_\delta \le  C \int_M | \nabla u_\delta |^2 \, dg_\delta.$$
The metrics $g_\delta$ converge uniformly to $g$, implying we may use a uniform constant $C$ in the Poincar{\'e} inequality for all $\delta$.  From this the first limit of the corollary follows.

A short computation yields
$$\int_M (u_\delta - \overline{u}_\delta)^2 \, dg_\delta = \vol(M,g_{\dd}) - \vol(M,g_{\dd}) \cdot \overline{u}_\delta^2.$$
The second limit of the corollary follows from this and the previous paragraph.
\end{proof}

Define $$w_\delta := u_\delta - \overline{u}_\delta.$$  From the fact that $u_\delta$ is an eigenfunction we obtain the following equation for $w_\delta$:
\begin{equation} \label{w}
-8 \Delta_\delta w_\delta + (\mathcal{R}_\delta - \lambda_\delta) w_\delta = (\lambda_\delta - \mathcal{R}_\delta) \overline{u}_\delta. 
\end{equation}
By the Sobolev inequality
$$\|w_{\dd}\|_{6} \leq C(\| \nabla u_{\dd}\|_{2}+ \| w_{\dd}\|_{2})$$
can be applied with a uniform constant $C$ for all the metrics $g_\delta$. By the above corollaries, the right hand side $\to 0$,
so
$$\| w_{\dd}\|_{6} \longrightarrow 0.$$
Using equation \ref{w}, it follows from the theory of elliptic partial differential equations \cite[Thm. 8.17]{GT01} that
$$\sup_M |w_\delta | \le C_2 \left\{ \left( \int_M w_\delta^6 \, dg_\delta \right)^{1/6} + \left( \int_M \overline{u}_\delta^3 |\lambda_\delta - \mathcal{R}_\delta |^3 \, dg_\delta \right)^{1/3} \right\} \longrightarrow 0.$$
It is again worth noting that a single constant $C_2$ is valid for all the metrics $g_\delta$ because they converge uniformly to $g$.  (Written in divergence form, the derivatives of the metric tensor do not appear in the coefficients of the local expression of the second order elliptic PDE satisfied by $w_\delta$.  So nothing is blowing up as $\delta \rightarrow 0$.)  We have shown that $u_\delta \rightarrow \overline{u}_{\dd} \rightarrow 1$ uniformly.

We may now define a new family of metrics $\U{g}_\delta$ on $M$ by the conformal change
$$\U{g}_\delta := u_\delta^4 \cdot g_\delta.$$
It is known that under such a conformal change, the scalar curvature of the new metric is given by the formula
$$R(\U{g}_\delta) = u_\delta^{-5} \left( -8 \Delta_\delta + R(g_\delta) \right) u_\delta$$
\cite{Schoen89}.  Since $u_\delta$ is an eigenfunction we have the equation
\begin{eqnarray*}
(-8 \Delta_\delta + R(g_\delta)) u_\delta & = & (\lambda_\delta - \mathcal{R}_\delta) u_\delta + R(g_\delta) u_\delta \\
& = & (R(g_\delta) - \mathcal{R}_\delta) u_\delta + \lambda_\delta u_\delta.
\end{eqnarray*}
By the construction of the function $\mathcal{R}_\delta$ it follows that $R(g_\delta) - \mathcal{R}_\delta$ is everywhere bounded below by -6.  Since $\lambda_\delta \rightarrow 0$ and $u_\delta \rightarrow 1$ uniformly, we may conclude that the scalar curvature of the metric $\U{g}_\delta$ is bounded below by a constant converging to $-6$ as $\delta \rightarrow 0$.

Finally let us homothetically expand each metric $\U{g}_\delta$ by a factor $1+ \eta(\delta)$, for $\eta(\delta) \rightarrow 0$, so that the resulting homothetically expanded Riemannian metrics have scalar curvature bounded below by $-6$.  Let us call this final family of metrics $G_\delta$.  Notice that the metrics $G_\delta$ are converging uniformly to $g$.  In particular, the volumes $\text{Vol}(M,G_\delta)$ converge to the original $\text{Vol}(M,g) = \text{Vol}(M_+, g_+) + \text{Vol}(M_-, g_-)$. This finishes the proof of proposition \ref{scalar prop}.

\section{Bounding $\l(g_{\dd})$}
This section gives an alternative to the argument in section \ref{scalar section}, by 
showing that we may approximate $(M,g)$ (constructed at the beginning of section
\ref{scalar section}) by metrics with $\l(g_{\dd})\geq -6$
rather than $R(G_{\dd})\geq -6$. As mentioned in section \ref{monotonicity section}, we are 
including these two alternatives since the monotonicity of $V_{R}(g)$ under
Ricci flow with surgery is much simpler to prove than the monotonicity of $V_{\l}(g)$. 

Let $(M,g)$ be the metric constructed from equation \ref{glue} and satisfying equation \ref{gluing condition}
from section 4. Let $g_{\dd}$ be the smooth metrics approximating  the metric $(M,g)$ from proposition \ref{lowerbound}, such
that in $\S\times [-\dd,\dd]$, $R(g_{\dd})=R_{\delta}\geq \mathcal{S}$, and in $M-(\S\times [-\dd,\dd])$,
$R_{\dd}\geq -6$, where $\mathcal{S}$ comes from proposition \ref{lowerbound}. Let $\l(g)$ denote the minimal eigenvalue of the operator
$-4\D_{g} +R(g)$.

\begin{theorem}
$\underset{\delta \to 0}{\lim\inf} \l(g_{\delta}) \geq -6$. 
\end{theorem}
\begin{proof}
If $\mathcal{S}\geq -6$, then the result follows directly from the fact that $\l(g_{\dd})\geq \inf R_{\dd}$, so
we may assume that $\mathcal{S} < -6$.  
Let $u_{\dd}$ be the eigenfunction of $-4\D_{g_{\dd}}+R_{\dd}$ with minimal 
eigenvalue $\l(g_{\dd})<0$, such that $u_{\dd}$ is strictly positive, and $\| u_{\dd}\|_{2}=1$.
We can use the Rayleigh quotient method to estimate $\l(g_{\dd})$ from below. 
We may estimate 
$$ \l(g_{\dd})=\int_{M} |\nabla_{\dd} u_{\dd}|^{2} + R_{\dd} u_{\dd}^{2} dg_{\dd} \geq 
\int_{M-(\S\times [-\dd,\dd])} R_{\dd} u_{\dd}^{2} dg_{\dd} +  \int_{\S\times [-\dd,\dd]} R_{\dd} u_{\dd}^{2} dg_{\dd} 
$$
$$\geq -6\int_{M-(\S\times [-\dd,\dd])} u_{\dd}^{2} dg_{\dd} + \mathcal{S} \int_{\S\times [-\dd,\dd]} u_{\dd}^{2} dg_{\dd} 
\geq -6  + (\mathcal{S}+6) \int_{\S\times [-\dd,\dd]} u_{\dd}^{2} dg_{\dd} $$
Thus, we need only show that $ \int_{\S\times [-\dd,\dd]} u_{\dd}^{2} dg_{\dd} \to 0$  as $\dd \to 0$. 
Applying Holder's inequality, we see that
$$\int_{\S\times[-\dd,\dd]} u_{\dd}^{2} dg_{\dd}  \leq \left(\int_{\S\times [-\dd,\dd]} dg_{\dd}\right)^{\frac23} \left(\int_{\S\times [-\dd,\dd]} u_{\dd}^{6} dg_{\dd}\right)^{\frac13} = \vol(\S\times [-\dd,\dd],g_{\dd})^{\frac23}  \|u_{\dd}\|_{6}^{2}.$$
The right hand term converges to $0$, since $\vol(\S\times [-\dd,\dd],g_{\dd}) \to 0$ as $\dd \to 0$,
and $\|u_{\dd}\|_{6}$ is bounded by the following.
Since the metrics $g_{\dd}$ approximate $g$ in  $C^{0}$ norm, there is a uniform Sobolev constant $C$ (see the
justification in lemma \ref{Sobolev}) such
that 
$$\|u_{\dd}\|_{6} \leq C( \| \nabla_{\dd} u_{\dd}\|_{2} + \| u_{\dd}\|_{2})\leq C\left((\l(g_{\dd})-\mathcal{S})^{\frac12}+1\right) \leq C((-\mathcal{S})^{\frac12}+1).$$
\end{proof}

Rescaling the metrics $g_{\dd}$ by factors converging to $1$, we may approximate $(M,g)$ by metrics $g'_{\dd}$ such that $\l(g'_{\dd})\geq -6$.

\section{Ricci flow with $C^{0}$ input}

Let $(M,g)$ be the smooth closed $3$-manifold equipped with the Lipschitz Riemannian metric $g$ from the beginning of section \ref{scalar section}.  Let $\{ g_\delta \}_{\delta>0}$ be a family of smooth Riemannian metrics on $M$ converging uniformly to $g$.  We do not yet place any other geometric assumptions on the approximates $g_\delta$.  The goal of this section is to establish facts about flowing the singular metric $g$ in the direction of its Ricci curvature.

The analytic tools of this section were proven by M. Simon \cite{MSimon02}, where the dual Ricci harmonic map heat flow is studied with nonsmooth initial data. To apply Simon's machinery it is necessary to fix a ``background metric''.  Concretely this means fixing a smooth Riemannian metric $h$ on $M$ so that the identity map $(M,h) \longrightarrow (M,g)$ is $K$-bilipschitz, where $K>1$ is a fixed constant depending only on the dimension.  For our purposes, it would suffice to pick a metric $g_\delta$ for sufficiently small $\delta$.  Let us choose such a metric and call it $h$. The following theorem follows from statements in section 5 of \cite{MSimon02}.

\begin{theorem} \label{Simon theorem} 
There exists a constant $T>0$ depending only on the background metric $h$, a sequence $\delta_j \rightarrow 0$, and a sequence of families of smooth metrics ${g(\delta_j,t)}, 0 \le t\le T, j \in \NN$ such that: \\
(1) The metric $g(0,0) =g$ and $g(\delta_j, 0) = g_{\delta_j}$.\\
(2) The family $\{ g(\delta_j,t) \}_{0\le t \le T}$ is conjugate via diffeomorphisms of $M$ to a solution to the Ricci flow equation with initial metric $g(\delta_j,0)$.  Similarly, for any $t_0 >0$ the family $\{ g(0,t) \}_{t_0 \le t \le T}$ is conjugate via diffeomorphisms of $M$ to a solution to the Ricci flow equation with initial metric $g(0, t_0)$. \\
(3) $g(0,t)$ converges uniformly to $g(0,0)=g$ as $t \rightarrow 0$. \\
(4) $g(\delta_j, t)$ converges in the $\mathcal{C}^\infty$ sense to $g(\dd_{j},0)$ as $t \rightarrow 0$.\\
(5) For $t>0$, $g(\delta_j, t)$ converges in the $\mathcal{C}^\infty$ sense to $g(0, t)$ as $j \rightarrow \infty$.\\
\end{theorem}
\begin{proof}

Properties (1), (2), and (3) above are stated in the proof  of theorem 5.2 of \cite{MSimon02}.  Property (4) is a restatement of  lemma 5.1 of \cite{MSimon02}.  Also in theorem 5.2 of \cite{MSimon02}, it is proven that there exists constants $\{c_{i}\}$ depending on $h$ such that
 $$\sup_M | \nabla^i g(\delta_j, t) |^2 \le c_i/ t^i,$$
where $\nabla^i$ denotes the $i^{\text{th}}$ covariant derivative with respect to the background metric $h$.  (By fixing a finite set of coordinate charts covering $M$, and choosing a precompact domain in each chart whose union still covers $M$, we can replace $\nabla^i$ with honest spatial derivatives in each coordinate chart.  This will change the constants $c_i$ depending on how we choose our cover.  This is possible because on each precompact domain the identity map from the metric $h$ to the Euclidean metric will be bilipschitz.)  Notice the constants of this inequality do not depend on $\delta_j$.  For any $t_0 > 0$ this independence yields uniform control on the spatial derivatives of the metrics $g(\delta_j, t)$ for all $j$ and $t \ge t_0$.  We may therefore apply the Arzela-Ascoli theorem repeatedly to deduce property (5) above.  (This is in fact how Simon produces the family $g(0,t)$.  See the first paragraph of the proof of theorem 5.2.)
\end{proof}

Recall that for a smooth Riemannian metric $h$, $\lambda(h)$ denotes the minimal eigenvalue of the operator $-4\D_{h} +R(h)$.

\begin{cor} \label{simon cor 1} 
If for some $\mathcal{S} \in \mathbb{R}$ and all $j\in \mathbb{N}$ the metric $g_{\delta_j}$ has scalar curvature bounded below by $\mathcal{S}$, then for all $t>0$ the Riemannian metric $g(0,t)$ has eigenvalue $\lambda(g(0,t)) \ge \mathcal{S}$.
\end{cor}
\begin{proof}
 Perelman has shown that the eigenvalue $\lambda$ is nondecreasing along the Ricci flow (see Cor.1.2 and section 2.2 of \cite{Per02}).  Therefore by properties (2) and (5) of theorem \ref{Simon theorem}, together with continuity of $\lambda$ in the $\mathcal{C}^\infty$ topology on metrics, we can deduce the statement of the corollary.
\end{proof}

\begin{cor} \label{simon cor 2}
If for all $j \in \mathbb{N}$ the metrics $g_{\delta_j}$ have scalar curvature bounded below by $-6$, then for all $t>0$ the Riemannian metric $g(0,t)$ has scalar curvature bounded below by $-6$.
\end{cor}
\begin{proof}
The minimum scalar curvature is nondecreasing along the Ricci flow.  Therefore the corollary follows from theorem \ref{Simon theorem}.
\end{proof}
{\bf Remark:} Theorem 6.6 of \cite{MSimon02} shows that a manifold with a Lipschitz metric and Alexandrov positive
curvature which is bounded above and below may be approximated by smooth manifolds with positive curvature. 
Simon mentions at the end of  theorem 6.6 that one may prove a similar statement for positive scalar curvature. 
But his theorem is weaker than corollary \ref{simon cor 2}, since we do not assume that the metric may be 
approximated by smooth metrics with uniformly bounded curvature. The advantage of dealing with scalar
curvature (or $\l$) is that it is monotonic under the Ricci flow, whereas almost positive curvature is not, which
accounts for our weaker hypotheses.

\section{Main inequality}

Before we state the next theorem, recall that $-6$ is the scalar curvature of a $3$-manifold with constant sectional curvature $-1$.  

\begin{theorem} \label{main inequality}
Let $(X,g)$ be a compact $3$-manifold with a smooth Riemannian metric $g$ such that the scalar curvature of $g$ is at least $-6$ at every point, and the boundary of $X$ is a minimal surface.  Then
$$\text{Vol}(X,g) \ge \frac{1}{2} \, V_{3}\|DX\|.$$
Equality holds if and only if $X$ has constant sectional curvature $-1$ and the boundary of $X$ is totally geodesic.
\end{theorem}
\begin{proof}
We present two slightly different proofs, one using scalar curvature and the other using the eigenvalue $\lambda(h)$ of the operator $-4 \Delta_h + R(h)$.  Let us discuss the scalar curvature proof first.

Metrically double the manifold $(X,g)$ across its boundary to form the closed manifold $M=DX$ equipped with the piecewise smooth Lipschitz continuous metric formed by two copies of $g$  (see equation \ref{glue} of section \ref{scalar section}.)  Abusing notation slightly, let us denote this Lipschitz Riemannian manifold by $(M,g)$.

By proposition \ref{scalar prop}, there is a family of smooth metrics $\{ g_\delta \}$ on $M$ converging uniformly to $g$ with scalar curvature bounded below by $-6$.  Now let $g(t)$ be the family of metrics produced by theorem \ref{Simon theorem} satisfying the properties: as $t \rightarrow 0$ the metrics $g(t)$ converge uniformly to $g$, and for any $t_0>0$ the family $\{ g_t \}_{t_0 \le t \le T}$ is diffeomorphism-conjugate to a solution to the Ricci flow with initial metric $g_{t_0}$.  By corollary \ref{simon cor 2}, for any $t>0$, $R_{min}(g_t)\geq -6$. We may assume 
$R_{min}(g_{t})<0$ for all $t$, otherwise we are done.  
By theorem \ref{monotonicity}, $V_{R}(g_{t}) \geq V_{3}\|M\|$.
Taking the initial time $t_0$ down to $0$, and using the fact that $\underset{t\to 0}{\lim \inf}\, V_{R}(g_{t}) \leq  \vol(M,g)$ by \ref{simon cor 2} yields the desired inequality
$$\text{Vol}(X,g) = \frac{1}{2} \, \text{Vol}(M,g) \ge V_{3}\|M\|.$$ In this argument, we may replace $V_{R}(g_{t})$ with $V_{\l}(g_{t})$, and use \ref{simon cor 1} instead.  

Let us now assume that $\text{Vol}(X,g) = \frac{1}{2} V_{3}\|M\|$.  Then it follows that the quantity $V_{R}(g_{t})$ is constant for all time.  This can only occur if all the metrics $\{ (M, g_t) \}_{t>0}$ have constant Ricci curvature \cite[Sec.2]{Ham99} (one may also use $V_{\l}(g_{t})$ is constant, and use \cite[Sec.2.3]{Per02} and \cite[Sec.8]{Per03}) .  In $3$ dimensions constant Ricci curvature implies constant sectional curvature.  By the uniform continuity of the metrics back to time zero, this implies that $(M,g)$ is isometric to a smooth Riemannian manifold of constant sectional curvature $-1$.  The gluing between the two copies of $(X,g)$ produces a smooth Riemannian metric along the boundary only if the boundary of $(X,g)$ is totally geodesic.  This completes the proof of the theorem.
\end{proof}

\section{Hyperbolic convex cores} \label{bonahon}

In this section we will use theorem \ref{main inequality} to study convex cores of hyperbolic $3$-manifolds.  More specifically we can quickly prove a conjecture of Bonahon stating that the volume of a compact hyperbolic $3$-manifold $M$ with non-empty convex boundary is at least half the simplicial volume of the doubled manifold $DM$.  In the incompressible boundary setting this was proved in \cite{Storm04}.  We will also prove the relative version of the inequality for cusped hyperbolic $3$-manifolds.

If $N$ is a complete (infinite volume) hyperbolic $3$-manifold without boundary, then the convex core $C_N \subseteq N$ is defined to be the minimal closed convex submanifold of $N$ for which inclusion is a homotopy equivalence.  The volume of the convex core $C_N$ is a measure of the geometric complexity of $N$.  When it has infinite volume then $N$ is said to be geometrically infinite.  Otherwise $N$ is geometrically finite.  This is the case of interest here.  

We need some standard geometric properties of the convex core which can be found in \cite{EM87}.  The distance function to the convex core defines a distance decreasing flow $\mathcal{F}_N$ on the complement of the convex core which simply flows each point to the nearest point in the convex core.  (In general this flow is only once differentiable, but this will suffice here.)  The boundary of the convex core is not in general smooth, but it is always totally geodesic outside some compact set.  The relationship between the topology of $N$ and $C_N$ is straightforward: if $N$ is geometrically finite then it is homeomorphic to the interior of $C_N$.

To begin we will prove the desired volume inequality in the compact setting. 

\begin{theorem}
Let $N$ be a complete hyperbolic $3$-manifold without boundary.  Assume the convex core $C_N$ of $N$ is compact.  It then satisfies the volume inequality
$$\text{Vol}(C_N) \ge \frac{1}{2} \, V_3 \, \| DC_N \|.$$
Moreover, if $2 \text{Vol}(C_N) = V_3 \, \| DC_N \|$ then $C_N$ is a compact manifold with totally geodesic boundary.
\end{theorem}

\begin{proof}
Consider the (possibly disconnected) compact embedded surface given by the boundary of the convex core.  By Dehn's lemma we can compress a finite number of homotopically trivial curves in $\partial C_N$ to obtain a $\pi_1$-injective embedded surface $S'$ bounding a submanifold $X'$. Throw away components $C$ of $X'$ such that $\|DC \|=0$ ({\it i.e.} $DC$ is a graph manifold) to get a submanifold $X''$ with boundary $S''$. It follows that $\|DX''\|=\|DC_{N}\|$, by \cite[Ch.6]{Th},\cite{Soma81}. Using the distance decreasing flow $\mathcal{F}_N$ to push wayward surfaces back into the convex core, one can produce an embedded minimal surface $S\subset C_N$ isotopic to $S''$ \cite{FHS83, HS}.  $S$ bounds a (possibly disconnected) $3$-manifold $X \subseteq C_N$ homeomorphic to $X''$.  We can apply theorem \ref{main inequality} to $X$ to obtain
$$\text{Vol}(X) \ge \frac{1}{2} \, V_3 \, \| DX \|.$$
If $X$ is not a proper submanifold of $C_N$ then $\partial C_N$ must be a minimal surface to begin with.  This is possible only if $\partial C_N$ is totally geodesic.  Therefore if $C_N$ is not totally geodesic then $X$ is a proper submanifold of the convex core, yielding the inequality
$$\text{Vol}(C_N) > \text{Vol}(X)\geq \frac12\, V_{3}\, \|DC_{N}\|.$$
\end{proof}

Proving the same inequality in the presence of cusps presents only technical difficulties.  The first is determining the proper statement of the theorem.  This we will attempt without introducing the language of pared $3$-manifolds.  (For a more detailed discussion of the cusped case see \cite{Storm04}.)  Let $N$ be a complete geometrically finite hyperbolic $3$-manifold without boundary.  Let $C_N$ be its convex core.  Assume that $C_N$ is not compact.  Then $DC_N$ is topologically a $3$-manifold with ends homeomorphic to a torus cross a ray.  (See \cite{Mo} for a precise description of the ends of $C_N$.)  Cut off the ends at an embedded $\pi_1$-injective torus to form a compact $3$-manifold with torus boundary.  Let us denote this compact manifold by $\overline{DC}_N$.  The Gromov norm rel boundary of a compact $3$-manifold with toroidal boundary is well defined \cite[Ch.6]{Th}.  Abusing notation slightly, let us denote this relative Gromov norm also by $\| \overline{DC}_N \|$.  (We will have no need of the non-relative Gromov norm for manifolds with boundary.  As one expects, $\| \overline{DC}_N \|$ equals half the Gromov norm of the closed manifold obtained by doubling $\overline{DC}_N$ along its boundary \cite[Ch.6]{Th},\cite{Soma81}.)  Using the product structure on $N- C_N$ \cite{EM87} and the precise manner which $C_N$ exits the cusps of $N$ \cite{Mo}, it can be shown that $\| \overline{DC}_N \|$ is an invariant of the relative homotopy type of $(N, \text{cusps of } N)$.  (This follows from \cite{Johannson79}.  See \cite{Storm04}.)

\begin{theorem} \label{noncompactbonahon}
Let $N$ be a complete geometrically finite hyperbolic $3$-manifold without boundary.  Assume that the convex core $C_N$ is not compact.  Let $\overline{DC}_N$ be the truncated double of $C_N$.  Then
$$\text{Vol}(C_N) \ge \frac{1}{2} \, V_3 \, \| \overline{DC}_N \|.$$
Moreover, if $ 2 \text{Vol}(C_N) = V_3 \, \| \overline{DC}_N \|$ then $C_N$ has totally geodesic boundary.
\end{theorem}
{\noindent}Note that if $N$ is geometrically infinite then by definition $\text{Vol}(C_N) = \infty$, so this is not a restrictive hypothesis.

\begin{proof}
Let us assume that the boundary of $C_N$ is not totally geodesic.  We will now describe a process for producing minimal surfaces in noncompact convex cores.  A similar procedure is described in \cite[Thm.4.4]{HRW99}.  

Take an exhausting sequence of submanifolds $\{ N_i \}$ of $N$ each of whose boundary is a collection of horospherical annuli, so that the boundary of $N_{i+1}$ is distance at least $D_{i+1}$ from the boundary of $N_i$, where $\{ D_i \}$ is a sequence we specify next.

By picking $D_i$ sufficiently large, the Riemannian metric on each $N_i$ can be altered near its boundary to produce a manifold $N'_i$ with the following properties:\\

{\noindent}1.  The metric is unchanged on $N_{i-1} \subset N'_i$.\\
{\noindent}2.  A collar of the boundary of $N'_i$ has a flat (i.e. zero curvature) product metric with totally geodesic annular boundary.\\
{\noindent}3.  The sectional curvature of $N'_i$ is between $-1- \varepsilon_i$ and zero, where $\varepsilon_i \rightarrow 0$. \\
{\noindent}4.  The volume of $(N'_i \cap C_N) - (N_{i-1} \cap C_N)$ goes to zero.\\
{\noindent}5.  The submanifold $N'_i - N_{i-1}$ has symmetry group given by isometries of the boundary annulus $\partial N_{i-1}$.\\

{\noindent}For the construction of such a metric see \cite{Leeb95,Bes00}.  (In proposition 2.3 of \cite{Leeb95} and in \cite{Bes00} an analogous construction is done in the finite volume setting.  Here, the geometry of the cusp is unchanged, and the topology is simpler.) 

By Dehn's lemma we can compress a finite number of homotopically trivial curves in $\partial C_N$ to obtain a $\pi_1$-injective embedded surface $S'' \subset C_N$.  Throw away components of $S''$ bounding $I$-bundle regions to get a surface $S'$. Double each manifold $N'_i$ along its boundary to produce a smooth Riemannian manifold $Z_i$.  The surface $S'\subset C_N$ doubles in $Z_i$ to a closed surface $\Sigma'_i \subset Z_i$.   There is a compact convex submanifold  $C_{Z_i}$ in each $Z_i$ consisting of the points originally coming from the two copies of the convex core of $N$ (the convexity follows from the convexity of $C_{N}$ and the symmetry assumption 5).  In each $Z_i$, there is a distance non-increasing flow $\mathcal{F}_{Z_i}$ defined on the complement of $C_{Z_i}$, which is the flow defined by the gradient of the distance function
to $C_{Z_{i}}$.  
Using the distance non-increasing flow $\mathcal{F}_{Z_i}$ to push wayward surfaces back into $C_{Z_i}$ (actually, just notice that the projection to $C_{Z_{i}}$ is an area decreasing map), one can produce an embedded minimal surface $\Sigma_i \subset C_{Z_i}$ isotopic to $\Sigma'_i$ \cite{FHS83, HS} (the surfaces cannot degenerate to an immersion, since we have thrown away all components of $S''$ bounding $I$-bundles).  Consider the intersection $S_i$ of the surface $\Sigma_i$ with the hyperbolic submanifold $N_{i-1}$ of $N'_i$.  $S_i$ is a stable minimal surface rel boundary inside $N_{i-1}$.  Therefore its curvature is less than or equal to $-1$ at interior points.

Consider now the sequence of surfaces $S_j \subset N_i$ for $j>i$.  By the Gauss-Bonnet theorem these surfaces have uniformly bounded area.  Also, the fact that they are stable minimal surfaces implies that for each $i$ their principal curvatures are uniformly bounded in $N_i$ \cite{Schoen83}.  They therefore form a precompact family.  Take $i$ to infinity and use a diagonalization argument to obtain a limit $S \subset C_N$ which is also a smooth minimal immersion \cite{HS} into the originial manifold $N$.  Since the original surface $S'$ was embedded, it follows that the limit surface $S$ is also embedded \cite{FHS83}. 

$S$ is isotopic to $S'$ (this is not hard to show, given the controlled geometry of $S$), so $S$ bounds a hyperbolic $3$-manifold $X \subset C_N$.  Recall we assumed that $\partial C_N$ is not totally geodesic.  It follows that $S$ is not equal to $\partial C_N$, so $\text{Vol}(C_N) > \text{Vol}(X)$.  Let $X_i \subset C_{Z_i}$ be the submanifold bounded by $\Sigma_i$.   Using the fact that the volume of $(N'_i \cap C_N) - (N_{i-1} \cap C_N)$ goes to zero, one can show that the volume of $X_i$ converges to twice the volume of $X$.  We may therefore pick an $i_0$ sufficiently large so that 
$$2 \text{Vol}(C_N) > ( 1 + \varepsilon_{i_0})^{3/2} \, \text{Vol}(X_{i_0}).$$
Here $-1-\varepsilon_{i_0}$ is the lower curvature bound for $N'_{i}$ from condition 3 above.  Define $c := \sqrt{1 + \varepsilon_{i_0}}$ and let $c X_{i_0}$ denote the space $X_{i_0}$ homothetically expanded by the factor $c$.  The scaled space $cX_{i_0}$ has curvature bounded below by $-1$.

Now apply theorem \ref{main inequality} to $cX_{i_0}$.  This yields the inequality
$$(1 + \varepsilon_{i_0})^{3/2} \, \text{Vol}(X_{i_0}) = \text{Vol}(cX_{i_0}) \ge \frac{1}{2} \, V_3 \, \| DX_{i_0} \|.$$
It remains only to prove the Gromov norm of $DX_{i_0}$ is twice that of $\overline{DC}_N$ (rel its boundary).  Starting with $\overline{DC}_N$ we can cut along a finite collection of embedded essential spheres, fill in the new boundary with balls, throwing away $I$-bundle components, and then double the resulting manifold along its toroidal boundary to obtain a manifold homeomorphic to $DX_{i_0}$.  Cutting along essential spheres or $S^{1}$ bundles does not alter the Gromov norm, and the Gromov norm adds when gluing two manifolds along essential boundary tori \cite[Ch.6]{Th},\cite{Soma81}.  Therefore the operation from $\overline{DC}_N$ to $DX_{i_0}$ doubles the Gromov norm.
\end{proof}

\section{Noncompact case}
We would like to extend theorem \ref{main inequality} to the case of non-compact manifolds with minimal
surface boundary. Unfortunately, we
can only extend this theorem to a particular case. 

\begin{theorem} \label{noncompact}
Let $\overline{N}$ be a compact manifold with interior $(N,g)$, a hyperbolic 3-manifold of finite volume. Let $\overline{\S}$ be
an embedded incompressible ($\pi_{1}$-injective) surface in $\overline{N}$. Then $$\vol(N,g)\geq \frac12 V_{3} \|D(N\split \S)\|.$$
\end{theorem}
\begin{proof}
We may proceed as in the proof of theorem \ref{noncompactbonahon}, and produce a sequence of manifolds $(Z_{i},g_{i})$ 
homeomorphic to $D(\overline{N})$ with $R(g_{i})\geq -6$, and such that $D(\overline{\S})$ is a minimal surface in $Z_{i}$. Then 
$\vol(Z_{i}\split D(\overline{\S})) \geq  \| D(N\split \S) \|$ by theorem \ref{main inequality}, and $\vol(Z_{i}) \to 2 \vol(N)$, so that the theorem follows. 
\end{proof}

\section{Lower bound on minimal volume orientable hyperbolic 3-manifold} \label{minvol}
The following theorem was proven in previous unpublished work of Agol
and Dunfield via a direct
smoothing argument, improving on the main result of \cite{Ag0}. Now we may apply proposition \ref{scalar prop} to reprove this. 

\begin{theorem}[Agol-Dunfield]  \label{bound}
Let $M$ be a hyperbolic 3-manifold with metric $\nu$. Let $\gamma
\subset M$ be a geodesic link in $M$ of length $l$ with an
embedded open tubular neighborhood $C$ of radius $R$, and with
complete hyperbolic metric $\tau$ on $M_\gamma=M-\g$. Then
$$\vol(M_\gamma,\tau) \leq (\coth 2R)^{3}
\left(\Vol(M,\nu)+\frac{\pi}{2} l \tanh R \tanh 2R\right)$$
\end{theorem}
\begin{proof}
We follow the method of \cite{Ag0}. We have a tubular neighborhood $C$ about
$\g$ of radius $R$. The mean curvature of $\partial C$
is $\kappa=\frac12(\coth R+\tanh R)=\coth 2R > 1$. Remove $C$, and  insert a horocusp with curvature
scaled so that the mean curvature on the boundary matches that of $\partial C$. 
Then we obtain a $C^{0}$ metric on $M_{\g}$ satisfying equation  \ref{glue}. Since a horotorus has mean curvature $1$, we 
rescale the hyperbolic metric on the horocusp to have sectional curvature $-\kappa^{2}$,
by scaling lengths in the standard hyperbolic metric by a factor of $\kappa^{-1}$, which makes
the horotori have mean curvature $\kappa$. The tube $C$ has volume $\pi l \sinh^{2}R$, 
and the area of $ \partial C=\pi l \sinh 2R$. The horocusp in the standard hyperbolic
metric has volume $\frac12$ the area of its boundary horotorus, so a horocusp with boundary
of mean curvature $\k$ will have volume $\frac{1}{2\k}$ times the area of its boundary. 
Thus, the volume of the horocusp will be $\frac{\pi l\sinh2R}{2\k}=\frac12 \pi l \sinh 2R \tanh 2R $.  We 
must rescale by $\k$ in order that the minimal scalar curvature is $-6$, obtaining a complete
$C^{0}$ metric $h$ on $M_{\g}$ satisfying the hypotheses of proposition \ref{scalar prop}. Then 
we have $$\vol(M_{\g},h) =\k^{3}\left(\vol(M,\nu)-\pi l \sinh^{2} R + \frac12\pi l \sinh 2R \tanh 2R \right)$$
$$ =  \k^{3}\left(\vol(M,\nu) + \pi l \sinh^{2} R\left(\frac{2 \cosh^{2}R}{\cosh 2R} -1\right)\right) =\k^{3}\left(\Vol(M,\nu)+\frac{\pi}{2} l \tanh R \tanh 2R\right)$$

To simplify things, we have slightly lied in the argument in the previous paragraph, since
Perelman's monotonicity result only applies to closed manifolds. To remedy this, we
perform high Dehn filling on $M_{\g}$, by removing $C$, and inserting
a tubular neighborhood of a geodesic in $\HH^{3}$ with the boundary slope along a 
sequence of meridians $\a_{i}$, $i\to \infty$. As the length of the meridian $\a_{i}$ goes to 
$\infty$, the  inserted tube will converge in the Gromov-Hausdorff limit to the  horotorus 
described in the previous paragraph, and we obtain the same estimate by taking
a limit of estimates obtained using these Dehn filled metrics (alternatively, one may use the 
doubling trick of theorem \ref{noncompactbonahon}).  
\end{proof}

The following result improves on proposition 5.4 of \cite{Pr3}. 
\begin{cor}
The minimal volume orientable hyperbolic 3-manifold $M$ has $\vol(M)\geq 0.67$.
\end{cor}
\begin{proof}
This follows as in the proof of proposition 5.4 \cite{Pr3}, using theorem \ref{bound} to replace theorem 5.2 \cite{Pr3}, and the fact that a 
minimal volume orientable hyperbolic 3-manifold has a tube of radius $\frac12\ln 3$ \cite{GMT}, as well as that the
minimal volume one cusped orientable hyperbolic 3-manifold has volume $\geq 2V_{3}=2.02988...$ \cite{CM}.
\end{proof}

{\bf Remark:} Theorem \ref{bound} allows us to
assail the problem of finding the minimal volume hyperbolic 3-manifold on various fronts, by either
classifying cusped orientable hyperbolic 3-manifolds with volume $<2.852$ (for example, by extending the arguments
of Cao and Meyerhoff \cite{CM}), or by classifying closed 3-manifolds with a minimal length geodesic with tube radius
$< 0.7$.

\section{Universal manifold pairings}
The paper \cite{FKNSWW05} defines the notion of a universal $(n+1)$-manifold pairing. This consists
of linear combinations of $(n+1)$ manifolds with marked boundary, together with a pairing on these manifolds
which have the same $n$-dimensional marked boundary which takes values in formal linear combinations
of closed $n+1$-manifolds.  A unitary TQFT may be thought of as a ``representation'' of such a pairing, assigning
a Hermitian vector space to each closed $n$-manifold, and number to each closed $(n+1)$-manifold, and satisfying
some extra axioms \cite{Atiyah90}.  
Conjecture 2.2 of \cite{FKNSWW05} states that the (2+1)-dimensional universal 3-dimensional pairing
faithfully detects linear combinations of 3-manifolds (there are no hermitian isotropic vectors).  
We provide some partial evidence for this conjecture in a special case. 

Let $\S$ be a closed oriented surface, and let  $\MM_{\S}$ denote the formal combinations of compact oriented 3-manifolds
$M$ with $\partial M=\S$ and coefficients in $\CC$. Define the pairing $(,):\MM_{\S}\times \MM_{\S}\to \MM$ by $(\sum_{i} a_{i} M_{i}, \sum_{j} b_{j} N_{j})= \sum_{i,j} a_{i}\overline{b}_{j} M_{i}\cup_{\S} \overline{N}_{j}$ (where $\MM=\MM_{\emptyset}$, and $\overline{N}_{j}$ denotes $N_{j}$
with the reversed orientation).  
The following conjecture generalizes in a natural way conjecture 2.2 of \cite{FKNSWW05}. 
\begin{conjecture}
If there exists $x_{1},...,x_{m}\in \MM_{\S}$ such that $\sum_{i} (x_{i},x_{i})=0\in \MM$, then $x_{1}=\cdots=x_{m}=0$. 
\end{conjecture}
We provide some evidence for this conjecture in a special case. 

\begin{cor} \label{utqft}
Suppose that each $x_{1},...,x_{m}\in \MM_{\S}$ is a formal linear combination of compact, acylindrical
3-manifolds, such that $\sum_{i} (x_{i},x_{i})=0\in \MM$. Let $x_{i}=\sum_{j\in I_{i}} x_{ij}M_{j}$, and 
suppose that $\sum_{i} |I_{i}|^{2}$ is minimal. Then for $j,k\in I_{i}$, 
$\|M_{j}\cup_{\S} \overline{M}_{k}\| = \| M_{j}\cup_{\S} \overline{M}_{j}\|$.
\end{cor}
\begin{proof}
Let $g_{j}$ denote the hyperbolic metric on $M_{j}$ with
totally geodesic boundary. Denote the manifold $M_{jk}=M_{j}\cup_{\S}\overline{M}_{k}$, $j,k\in I_{i}$. Let $\S_{jk}\subset M_{jk}$ be a least area representative of $\S\subset M_{jk}$ in the canonical hyperbolic metric on $M_{jk}$. Then $\|M_{jk}\| \geq
\frac12(\|M_{jj}\|+\|M_{kk}\|)$ by theorem \ref{main inequality}. Since $\sum_{i,j} |x_{ij}|^{2} M_{jj} \neq 0$, 
in order for $\sum_{i} (x_{i},x_{i})=0$, we must have $M_{aa}=M_{bc}$, for some $b\neq c$, $b,c \in I_{i}$, where $\| M_{aa}\|$ is minimal over $\|M_{jj}\|, j\in \cup_{i} I_{i}$.
Thus, $\|M_{aa}\|=\|M_{bc}\| \geq \frac12(\|M_{bb}\|+\| M_{cc}\|)\geq \|M_{aa}\|$ (the second inequality follows from the assumed
minimality of $\|M_{aa}\|$), and thus $\| M_{aa}\|=  \| M_{bb}\|= \|M_{cc}\|=\|M_{bc}\|$. Since we have equality,
we must have $\S_{bc}$ is
totally geodesic in $M_{bc}$ by Thm. \ref{main inequality}, and thus $g_{b}$ and $g_{c}$ induce the same metric up to isotopy on $\S$. Take the maximal sub partition $J=\{J_{1},...,J_{n}\}$ of $I=\{I_{1},...,I_{m}\}$  such that  for all $i,j\in J_{l}$,  $\|M_{ij}\|=\|M_{aa}\|$. This partition is a  subset of the partition by 
the isotopy class of the metric on $\S=\partial M_{j}$ induced by $g_{j}$, and therefore it is easy to see that $\cup J_{l}$ includes
all $j$ such that $\|M_{jj}\|=\|M_{aa}\|$. Then for $i$ not in $\cup J$, $i\neq j$, $\| M_{ij}\| > \|M_{aa}\|$, and thus all the cancellation among the
terms of $\sum_{i} (x_{i},x_{i})$ with $M_{aa}$ must occur among the terms of the partition $J$. Associated to $J_{l}$, there is a vector $y_{l}$, where $y_{l}= \sum_{j\in J_{l}} x_{ij} M_{j}$, where $J_{l}\subset I_{i}$. Then $\sum_{l} (y_{l},y_{l}) = 0$, and $\sum_{l} |J_{l}|^{2}\leq \sum_{i} |I_{i}|^{2}$. By the minimality hypothesis of $I$, we have $J=I$, and the result follows. 
\end{proof}

This reduces  this special case of conjecture 2.2 of \cite{FKNSWW05} to a geometric question about how many distinct ways a hyperbolic
manifold may be obtained by doubling a manifold with geodesic boundary. It seems promising
that this question may yield to geometric techniques.

\section{Conclusion}
The results in this paper give rise to many interesting questions. A natural question is whether one
may prove the main result without using Ricci flow.
\begin{conjecture}
If a finite volume hyperbolic 3-manifold with minimal surface boundary is  locally minimal among such
manifolds, then the boundary is totally geodesic.
\end{conjecture}
If this conjecture were true, then one should be able
to reprove theorem \ref{main inequality} by deforming a hyperbolic manifold with minimal surface
boundary to have geodesic boundary, while decreasing volume. 

Another possible approach would be to try to use the methods of the natural map \cite{BCG}. If $(M,g)$ is
a Riemannian metric, let $h(M,g)$ denote the volume entropy of $g$. 
\begin{conjecture}
If $(M,g)$ is a closed Riemannian 3-manifold with $R(g)\geq -6$, then $h(M,g)\leq h(\HH^{3})$. 
\end{conjecture}
This would enable one to reprove the main theorem using the techniques of the natural map. 

One nagging point we were unable to resolve is to prove the main theorem in the case
that the manifold is non-compact. 
\begin{conjecture}
Let $(M,g)$ be a complete finite volume 3-manifold with minimal surface boundary and scalar curvature $R(g)\geq -6$. 
Then $\vol(M,g)\geq \frac12 \|DM\|$, with equality iff $M$ has geodesic boundary. 
\end{conjecture}
This conjecture might be useful in an attempt to prove conjecture 2.2 of \cite{FKNSWW05}, by extending
the arguments of corollary \ref{utqft}.


\def\cprime{$'$} \def\cprime{$'$} \def\cprime{$'$}
\providecommand{\bysame}{\leavevmode\hbox to3em{\hrulefill}\thinspace}
\providecommand{\href}[2]{#2}

%
%
%
%
%


\vfill
\pagebreak

\appendix

\section{Volume change under drilling: theory vs.~experiment}
\centerline{\textsc{by Nathan M. Dunfield}}

\vspace{0.3cm}

\newcommand{\Vfill}{V_{\mathrm{fill}}}
\newcommand{\Vdrill}{V_{\mathrm{drill}}}
\newcommand{\Vest}{V_{\mathrm{est}}}
\newcommand{\Corig}{C_O}
\newcommand{\Cperel}{C_{P}}

Let $M$ be a closed hyperbolic 3-manifold, and let $\gamma$ be a simple
geodesic in $M$.  Consider the complement $M_\gamma = M \setminus \gamma$.  This
appendix focuses on the following question: how are the volumes of $M$
and $M_\gamma$ related?  In general, Thurston showed that $\vol(M) <
\vol(M_{\gamma})$, so the goal here is to control the increase in volume
when we drill out $\gamma$.  To be able to do this, we need more
information about the geometry of $M$, so assume in addition we know the
length $L$ of $\gamma$, and the radius $R$ of a maximal embedded open
tube about $\gamma$.  (Equivalently, $R$ is half the minimal distance
between two distinct lifts of $\gamma$ to the universal cover of $M$.)
To simplify the notation, set $\Vfill = \vol(M)$ and $\Vdrill =
\vol(M_\gamma)$.  The volume estimates I will discuss involve the
term:
\[
B = \Vfill + \pi L \sinh^2(R) \sech(2 R).
\]
In \cite{Ag0}, Agol showed that
\[
\Vdrill \leq  \Corig B \quad \mbox{where} \quad  \Corig = \big(\coth(R) \coth(2 R)\big)^{3/2}.
\]
Using Perelman's work, Agol and I improved this to
\begin{equation}\label{newest}
\Vdrill \leq \Cperel B  \quad \mbox{where} \quad \Cperel = \coth^3( 2 R) .
\end{equation}
See Theorem~\ref{bound} in the body of this paper for a proof.  The
following plot shows the ratio $\Corig/\Cperel$ of the multiplicative
factors as a function of $R$:

\vspace{0.5cm} 

\begin{center}
\begin{overpic}[tics=10, scale=0.8]{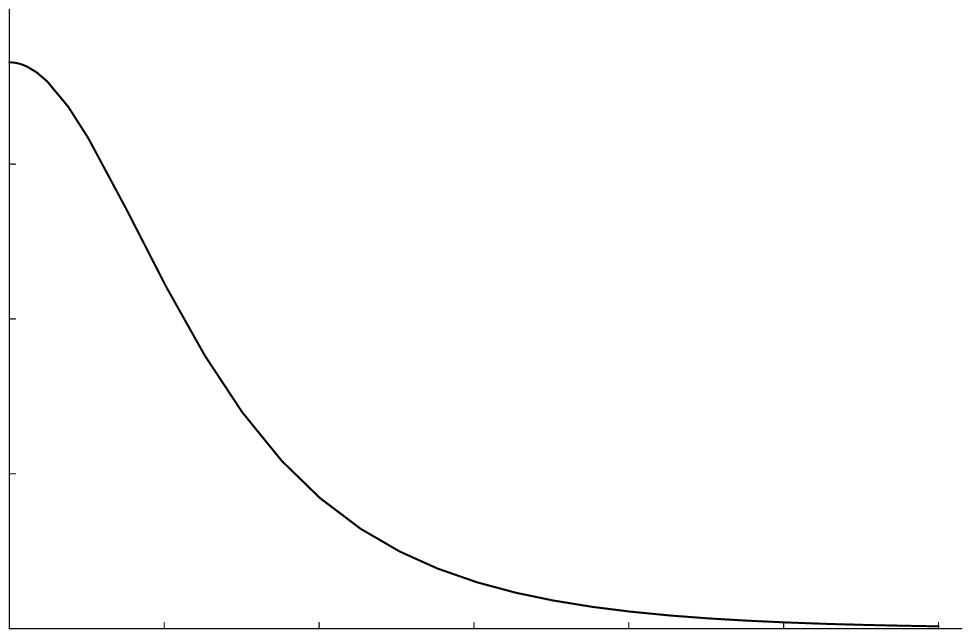}
\put(0,-5){$0$}
\put(13,-5){$0.5$}
\put(31.5,-5){$1$}
\put(45,-5){$1.5$}
\put(63,-5){$2$}
\put(78,-5){$2.5$}
\put(96,-5){$3$}
\put(-24, 30){${\displaystyle \frac{\Corig}{\Cperel}}$}
\put(35,-14){Tube radius $R$}
\put(0,-5){0}
\put(-4,-1){1}
\put(-8,15){1.5}
\put(-4,31){2}
\put(-8,47){2.5}
\end{overpic}

\vspace{1.5cm}

\end{center}

\noindent
As you can see, the new estimate is markedly better than the old one for
small $R$, but the two estimates are asymtopically the same for large $R$.

The purpose of this appendix is to compare the new drilling estimate
with data from more than $25{,}000$ pairs $(M, \gamma)$.  In all cases,
the data satisfies (\ref{newest}).  Moreover, the estimate in
(\ref{newest}) is fairly good for $R > 1$.  Since (\ref{newest}) rests
on Perelman's announced proof of geometrization \cite{Per02}, these
results can be viewed as an \emph{experimental} check on the
correctness of his proof.  There are two obvious caveats here: only a
finite number of manifolds were examined, and only part of Perelman's
work is needed (namely Theorem~\ref{monotonicity} in the body of this
paper).  There is a more subtle issue as well: it is entirely possible
that for these examples, when one applies the proof of
Theorem~\ref{monotonicity} the Ricci flow never goes through any
surgeries.  That is, these examples could well avoid one of the
least-understood parts of Perelman's approach.

After briefly discussing the sample manifolds, I will give various
graphical comparisons between (\ref{newest}) and the data.  I will
then turn to an alternate way of explaining the data, namely
hyperbolic Dehn surgery.  This approach more closely models
the data than (\ref{newest}), but does not provide universal bounds as
(\ref{newest}) does.  

\subsection{The sample}

The sample consists of $25{,}709$ pairs $(M, \gamma)$, where the
manifolds $M$ were drawn from the Hodgson-Weeks census \cite{SnapPea}.
The geometry of selected geodesics were computed using Oliver
Goodman's program \emph{Tube}, which is distributed with his program
\emph{Snap} \cite{Snap}.  The output of the program \emph{Tube} is not
completely rigorous, but should be essentially correct for nearly all
of the examples.  Most of the geodesics are not too long $(L \leq 2.5)$,
and have large tube radii $(R \geq 0.4)$.  Concentrating on large tube
radii was forced because \emph{Tube} had difficulty with the longer
geodesics.  However, as you will see, the estimate (\ref{newest}) is
really poor even for modestly small $R$, so this omission is of no
real importance.  While most of the geodesics are fairly short, to
concentrate on the regime where the two estimates differ noticeably, I
did require the geodesics to have length $L \geq 0.3$.  In more than
$24{,}000$ of the cases, the manifold $M_\gamma$ has a triangulation
with $\leq 7$ tetrahedra.  Partially as a result of having such simple
complements, I expect that, in nearly all cases, one can do hyperbolic
Dehn filling from the drilled manifold to the closed one.

\subsection{Experimental results}

\begin{figure}[htb]
\begin{center}
  \begin{overpic}[width=4in]{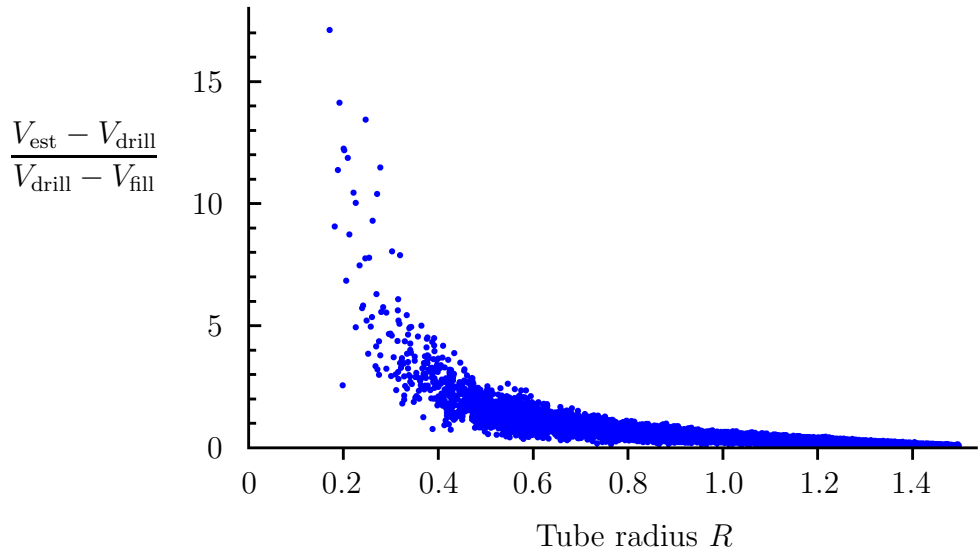}
    \put(1.5,-3){$0$}
    \put(12,-3){$0.2$}
    \put(24.7,-3){$0.4$}
    \put(37,-3){$0.6$}
    \put(49.5,-3){$0.8$}
    \put(62,-3){$1.0$}
    \put(74.5,-3){$1.2$}
    \put(86.5,-3){$1.4$}
    \put(-3, 1.5){$0$}
    \put(-3, 17.5){$5$}
    \put(-5, 33.5){$10$}
    \put(-5, 49.5){$15$}
    \put(-29, 40){${\displaystyle \frac{\Vest - \Vdrill}{\Vdrill - \Vfill}}$}
    \put(40, -10){Tube radius $R$}
  \end{overpic}
\end{center}
\vspace{1cm}
\caption{
  This plot shows the amount that (\ref{newest}) overshoots the actual
  volume $\Vdrill$, compared via (\ref{eq-first-comp}).  Three data
  points are omitted; in each case $R < 0.15$ and the vertical
  coordinate is $> 40$.  }  \label{fig-comp-1}
\end{figure}

\begin{figure}
  \vspace{0.2cm}
  \begin{center}
    \begin{minipage}{2.5in}
      \begin{overpic}[width=2.5in]{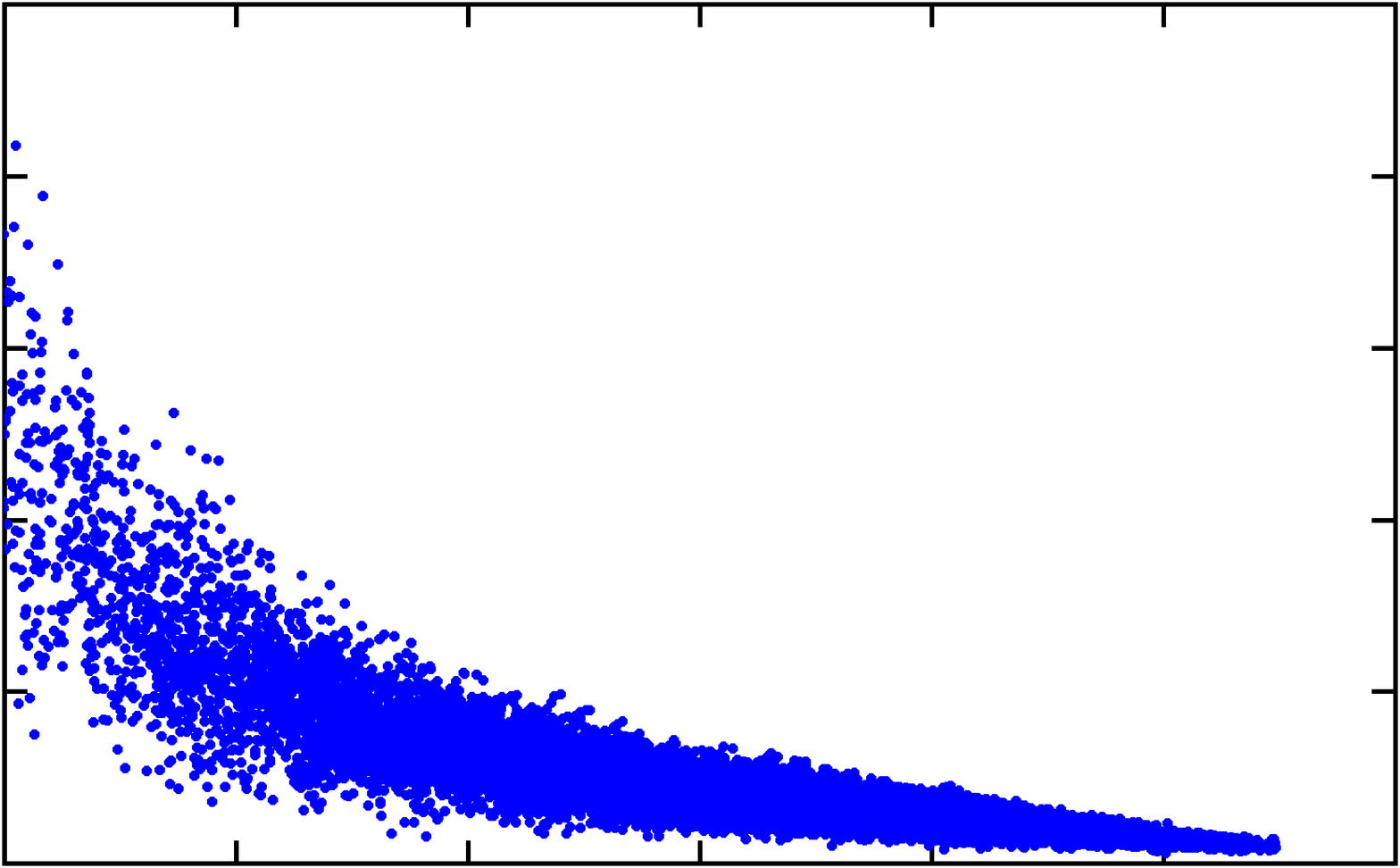}
        \put(13,-6){$0.6$}
        \put(29.5,-6){$0.8$}
        \put(46,-6){$1.0$}
        \put(63,-6){$1.2$}
        \put(79,-6){$1.4$}
        \put(-5, -1){$0$}
        \put(-5, 11){$1$}
        \put(-5, 22){$2$}
        \put(-5, 35){$3$}
        \put(-5, 48){$4$}
        \put(-5, 59.5){$5$}
      \end{overpic}
    \end{minipage}
    \hspace{1.5cm}
    \begin{minipage}{2.5in}
      \begin{overpic}[width=2.5in]{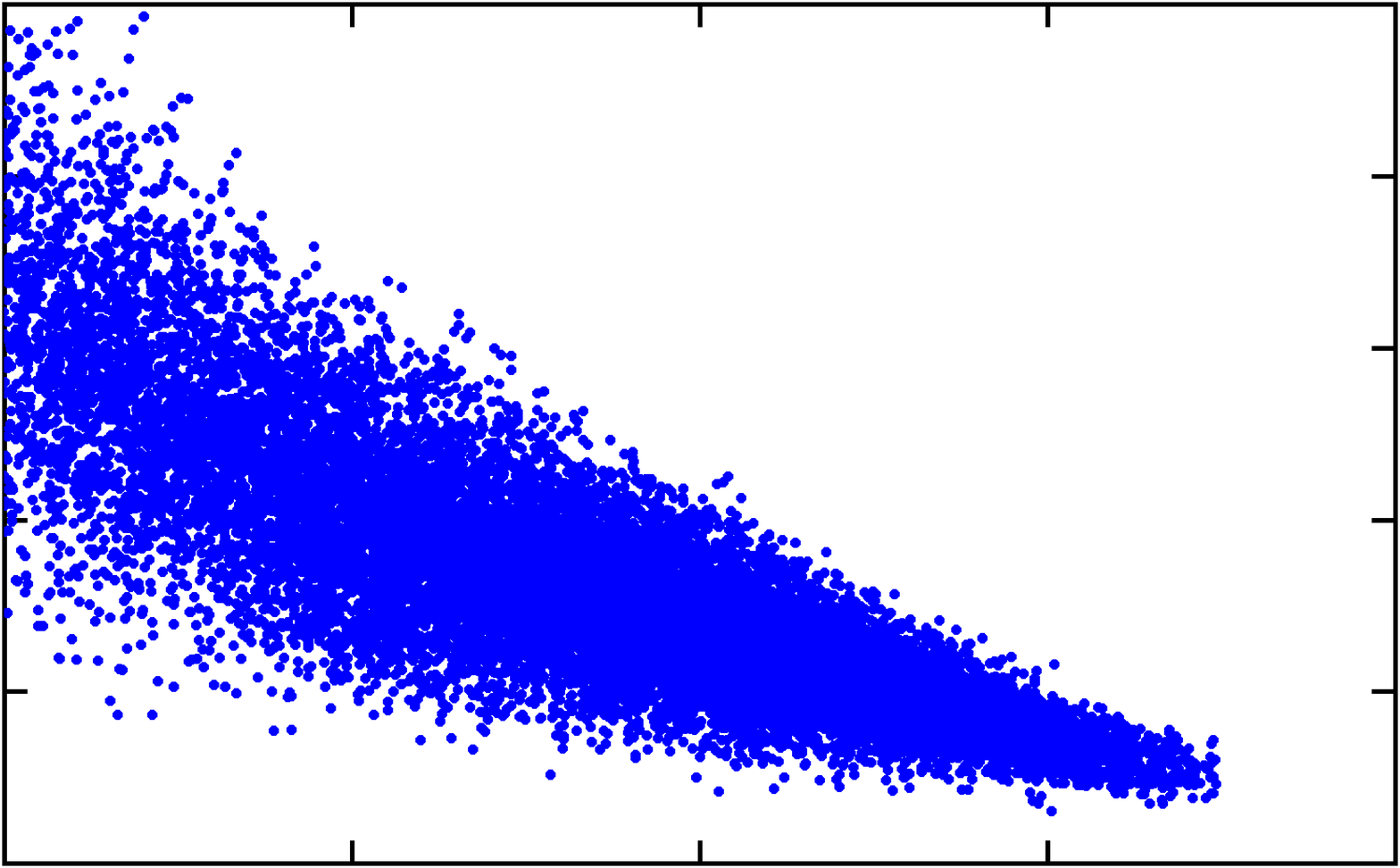}
        \put(21.5,-6){$1.0$}
        \put(46,-6){$1.2$}
        \put(71,-6){$1.4$}
        \put(-5, -1){$0$}
        \put(-10, 11){$0.2$}
        \put(-10, 22){$0.4$}
        \put(-10, 35){$0.6$}
        \put(-10, 48){$0.8$}
        \put(-5, 59.5){$1$}
      \end{overpic}
    \end{minipage}
    \vspace{0.3cm}
  \end{center}
  \caption{
    The data just for larger tube radii $R$.  The horizontal axis
    remains $R$, and the vertical axis the quantity
    (\ref{eq-first-comp}).  }\label{fig-comp1-close}
\end{figure}

\begin{figure}
  \vspace{0.18cm}
  \begin{center}
      \begin{overpic}[width=2.5in]{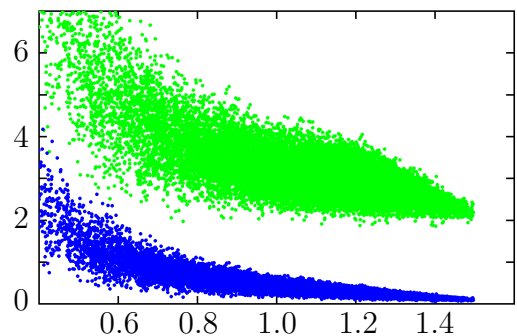}
        \put(13,-6){$0.6$}
        \put(29.5,-6){$0.8$}
        \put(46,-6){$1.0$}
        \put(63,-6){$1.2$}
        \put(79,-6){$1.4$}
        \put(-5, -1){$0$}
        \put(-5, 16){$2$}
        \put(-5, 33){$4$}
        \put(-5, 51){$6$}
      \end{overpic}
\end{center}
\vspace{0.3cm}
\caption{
  The top band is the estimate using the factor $C_O$,
  compared to the improved version using $C_P$ (bottom band).  As
  before, the horizontal axis remains the tube radius $R$, and the
  vertical axis the quantity (\ref{eq-first-comp}).  }

\label{fig-old-vs-new}
\end{figure}

Let $(M, \gamma)$ be as usual.  As per (\ref{newest}) set $\Vest = \Cperel B$.
Then we have
\[
        \Vfill < \Vdrill \leq \Vest.
\]
To understand how good a bound $\Vest$ is, we consider the overshoot error
$\Vest - \Vdrill$ as a proportion of the actual increase in volume
$\Vdrill - \Vfill$.   That is, in Figure~\ref{fig-comp-1} we plot
\begin{equation}\label{eq-first-comp}
  \frac{\Vest - \Vdrill}{\Vdrill - \Vfill}
\end{equation}  
against the tube radius $R$.  (One could just compare the error with
$\Vdrill$, but that provides less information when the tube radius is
large, in which case $\Vdrill$ is only a tad higher than $\Vfill$.)
Note that the vertical scale is quite large --- for the small $R$, the
error quantity (\ref{eq-first-comp}) is more than 15.  Thus for small
$R$, the estimate is quite poor.  This is not so surprising --- it
makes intuitive sense that the metric constructed in the proof of
(\ref{newest}) is far (in any reasonable sense) from the
hyperbolic metric on $M_\gamma$ when $R$ is small.  Thus one must
expect a larger drop in the normalized volume $V_R(M_\gamma)$ during the
Ricci flow with surgery, making the gap in Theorem~\ref{monotonicity}
larger.  Now, let's look at closeups for
larger $R$, where the estimate is markedly better; these are shown in
Figure~\ref{fig-comp1-close}.
It's also worth comparing Agol's original estimate compared to the Perelman
enhanced one; this is done in Figure~\ref{fig-old-vs-new}.
\FloatBarrier

\begin{figure}[h]
\begin{center}
  \begin{overpic}[width=4in]{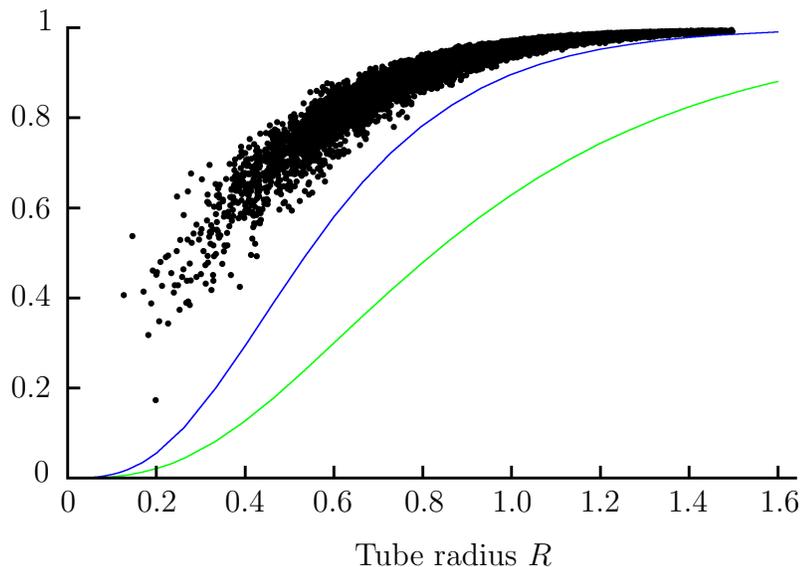}
    \put(1.5,-3){$0$}
    \put(11.5,-3){$0.2$}
    \put(23,-3){$0.4$}
    \put(34.5,-3){$0.6$}
    \put(46.5,-3){$0.8$}
    \put(58,-3){$1.0$}
    \put(69.5,-3){$1.2$}
    \put(81,-3){$1.4$}
    \put(93,-3){$1.6$}
    \put(-2, 1){$0$}
    \put(-5, 12){$0.2$}
    \put(-5, 24){$0.4$}
    \put(-5, 35.5){$0.6$}
    \put(-5, 47.5){$0.8$}
    \put(-1.5, 60){$1$}
    \put(40, -10){Tube radius $R$}
  \end{overpic}
\end{center}
\vspace{1cm}
\caption{
  Here, the vertical axis for the data is $B/\Vdrill$.  For
  (\ref{newest}) to hold, the data points must lie above
  $1/\Cperel(R)$ which is the upper curve.  The light lower curve is
  $1/\Corig(R)$.
  }\label{alt-comp-1}
\end{figure}

\begin{figure}[h]
  \begin{center}
    \begin{minipage}{2.5in}
      \begin{overpic}[width=2.5in]{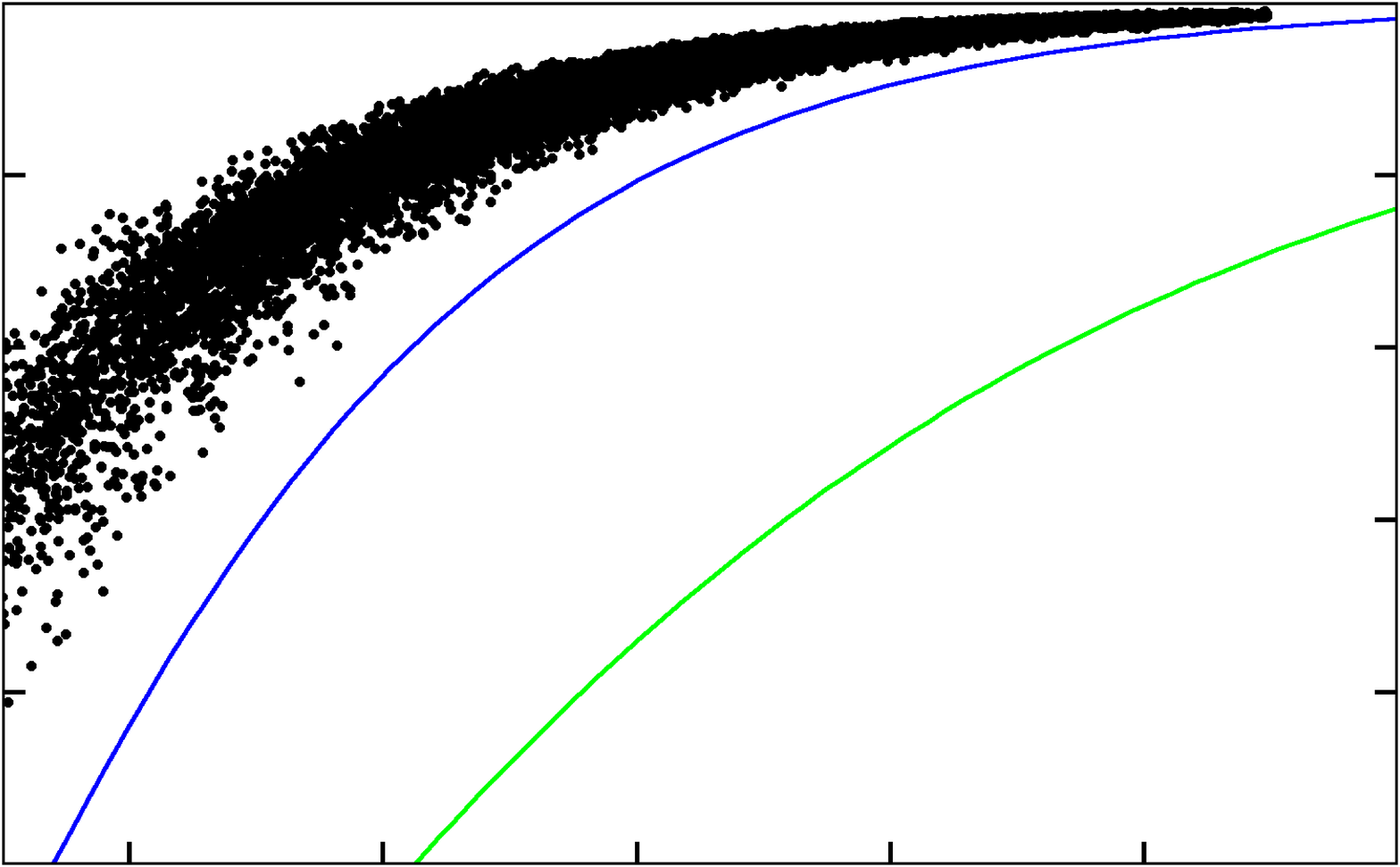}
        \put(5,-6){$0.6$}
        \put(23.5,-6){$0.8$}
        \put(41,-6){$1.0$}
        \put(60,-6){$1.2$}
        \put(77.5,-6){$1.4$}
        \put(-10, 11){$0.6$}
        \put(-10, 22){$0.7$}
        \put(-10, 35){$0.8$}
        \put(-10, 48){$0.9$}
        \put(-5, 59.5){$1$}
      \end{overpic}
    \end{minipage}
    \hspace{1.5cm}
    \begin{minipage}{2.5in}
      \begin{overpic}[width=2.5in]{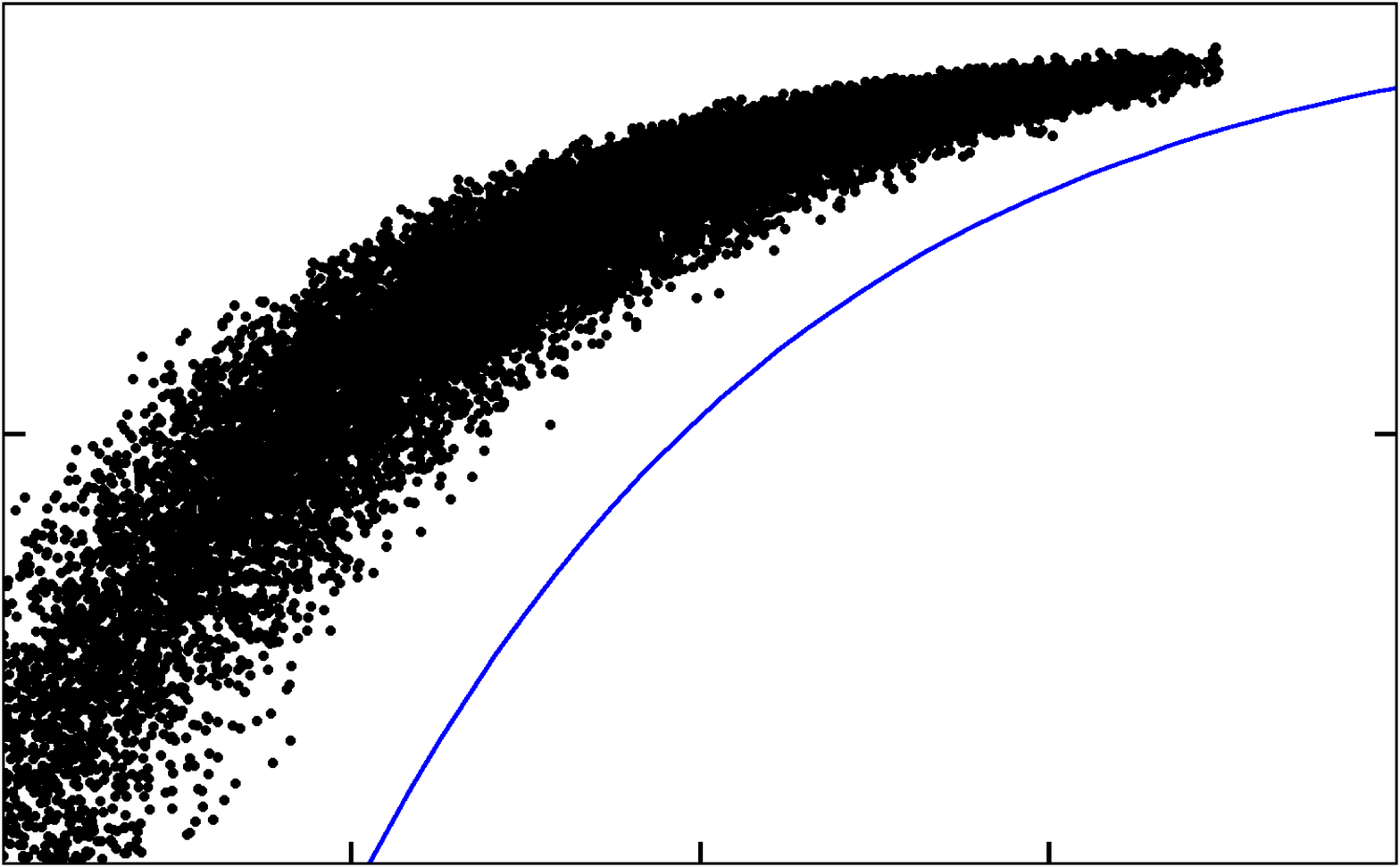}
        \put(21,-6){$1.0$}
        \put(46,-6){$1.2$}
        \put(71,-6){$1.4$}
        \put(-10, -1){$0.9$}
        \put(-13, 30){$0.95$}
        \put(-5, 59.5){$1$}
      \end{overpic}
    \end{minipage}
  \end{center}
  \vspace{0.3cm}
  \caption{Closeups of the data in Figure~\ref{alt-comp-1}. }\label{alt-comp-2}
\end{figure}

Another way to compare the data with the estimate is this.  Both of
the bounds have the same form $\Vdrill \leq C(R) B$.  It is natural to
ask: experimentally, what is the optimal form of the function $C(R)$ in
such an estimate?  Equivalently, we seek the function that just barely
satisfies:
\[ 
     \frac{1}{C(R)} \leq  \frac{B}{\Vdrill}
\]   
Figures~\ref{alt-comp-1}--\ref{alt-comp-2} show the left-hand side of
the above, plotted against $1/\Cperel$ and $1/\Corig$. 

One intriguing thing here is that $B/\Vdrill$ is always $\leq 1$.  In
other words,
\[
\Vfill + \pi L \sinh^2(R) \sech(2 R) \leq \Vdrill.
\]
This is not predicted by Theorem~\ref{bound}, and I don't see any
reason why it should hold true in general.  However, I also know of no
counterexamples.

\subsection{Hyperbolic Dehn surgery and volume change}

If the hyperbolic metric on $M$ can be deformed to the hyperbolic
metric on $M_\gamma$ through a series of cone manifolds with cone
locus $\gamma$, we say that $M$ is obtained from $M_\gamma$ by
\emph{hyperbolic Dehn filling}.  In this case, Schl\"afli's formula
gives that
\begin{equation}\label{Schalfi}
\Delta V = \Vdrill - \Vfill = \frac{1}{2} \int_0^{2 \pi} L(\theta) \, d \theta,
\end{equation}
where $L(\theta)$ is the length of $\gamma$ in the metric with cone angle $\theta$.

If we consider different hyperbolic Dehn fillings on a fixed cusp
manifold, Neumann and Zagier \cite{NZ} showed that $\Delta V \sim \pi L /2$ as the
length $L$ of the core geodesic goes to $0$.  For a particular Dehn
filling, Bridgeman \cite{Bri} observed that as long as $L(\theta)$ is an increasing
function of $\theta$, then (\ref{Schalfi}) implies that the volume
increase is bounded by twice the asymptotic estimate:
\begin{equation}\label{Bridgeman}
\Delta V \leq \pi L
\end{equation} 
While it appears that Bridgeman's inequality does not hold in general
\cite[Sec.~3]{Ag0}, Hogson and Kerckhoff showed it does hold when $L \leq
0.16$, provided also that the tube radius $R \geq 0.66$ \cite[Sec.~6]{HK}.  In the same
paper, Hogson and Kerckhoff gave a very precise version of the
asymptotics of $\Delta V \sim \pi L /2$ with similar hypotheses.  Given these
results, it makes sense to plot $\Delta V / (\pi L)$ as a function of $L$:
  \begin{center}
    \vspace{0.3cm}
      \begin{overpic}[width=4.5in]{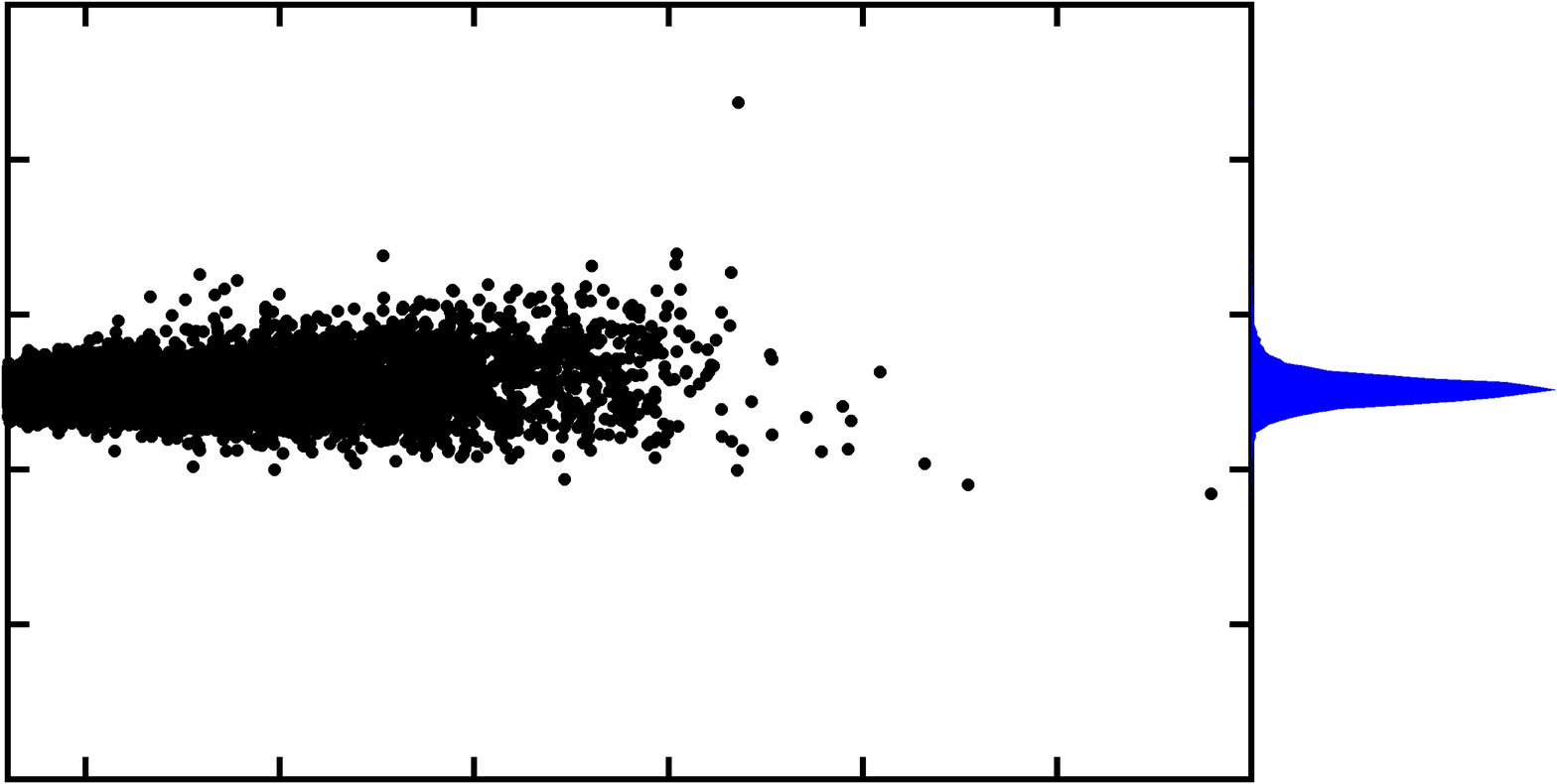}
        \put(5,-3){$0.5$}
        \put(18.3,-3){$1$}
        \put(29,-3){$1.5$}
        \put(42,-3){$2$}
        \put(52.7,-3){$2.5$}
        \put(66,-3){$3$}
        \put(-4, 9.5){$0.2$}
        \put(-4, 19){$0.4$}
        \put(-4, 28.5){$0.6$}
        \put(-4, 38.5){$0.8$}
        \put(-1, 48){$1$}
        \put(30, -10){Length $L$ of $\gamma$}
        \put(-15, 25){${\displaystyle \frac{\Delta V}{\pi L}}$}
      \end{overpic}
 \end{center}
 \vspace{1.5cm} As you can see, the result is the simplest picture we
 have seen yet --- the distribution of $\Delta V / (\pi L)$ is nearly
 independent of $L$, with some broadening as we increase $L$. The
 histogram attached to the right side of the plot gives the
 distribution of $\Delta V / (\pi L)$, which is a bell-curve shape with
 mean $0.5034$ (essentially the $1/2$ expected!)~and standard
 deviation $\sigma = 0.022$.  Indeed, the distribution is very nearly the
 same as the normal distribution with mean $0.502$ and $\sigma = 0.017$;
 the main difference is the top tail is a little thicker than the
 normal distribution.  In all cases, Bridgeman's inequality
 (\ref{Bridgeman}) holds.
 
 While the picture coming from hyperbolic Dehn filling provides the
 most accurate fit for our data, I should emphasize that it is not
 known whether one can always obtain $M$ from $M_\gamma$ via
 hyperbolic Dehn filling.  Even when this is known, one can only
 estimate $\Delta V$ if one makes additional assumptions about the behavior
 of $L(\theta)$.  The estimate (\ref{newest}) has the advantage that it
 applies without any additional hypotheses, which is important in
 applications such as bounding the volume of the smallest closed
 hyperbolic 3-manifold.

\end{document}